\documentclass[11pt]{article}
\usepackage{amssymb}
\usepackage{graphicx}
\usepackage{color}

\usepackage{amsmath}

\usepackage{amsfonts}

\usepackage{amssymb}
\usepackage{graphics}
 \usepackage{bbm}

\usepackage{titlesec}
\usepackage[titletoc,toc,title]{appendix}

\titleformat{\section}{\large\bfseries}{\thesection.}{0.3em}{}
\titleformat{\subsection}{\bfseries}{\thesubsection.}{0.3em}{}
\titleformat{\subsubsection}{\small\bfseries}{\thesubsubsection.}{0.3em}{}

 \allowdisplaybreaks[2]

\usepackage{latexsym}

\renewcommand{\skew}{\mathop{\rm skew}}
 \DeclareMathOperator{\sym}{sym}
\DeclareMathOperator{\tr}{tr} 
 \DeclareMathOperator{\dev}{dev}
\DeclareMathOperator{\Curl}{Curl\hskip.04truecm}

\DeclareMathOperator{\Div}{Div} 

\DeclareMathOperator{\sL}{\mathfrak{sl}}
\DeclareMathOperator{\so}{\mathfrak{so}}

\newcommand{\Chi}{\raisebox{0.5ex}{\mbox{{\Large $\chi$}}}}

\newcommand{\yieldlimit}{{\sigma}_{\mathrm{y}}}

\newcommand{\yieldzero}{{\sigma}_{0}}

\newcommand{\C}{\mathbb{C}}

\newcommand{\vectgam}{\raisebox{0.2ex}{\underline{$\gamma$}}}
\newcommand{\vecteta}{\raisebox{0.2ex}{$\underline{\eta}$}}

\newcommand{\BBR}{\mbox{$\mathbb{R}$}}

\newcommand{\SFH}{\mbox{$\mathsf{H}$}}

\newcommand{\SFV}{\mbox{$\mathsf{V}$}}
\newcommand{\SFW}{\mbox{$\mathsf{W}$}}
\newcommand{\SFP}{\mbox{$\mathsf{P}$}}
\newcommand{\SFZ}{\mbox{$\mathsf{Z}$}}
\newcommand{\bzero}{\mbox{$\bf 0$}}
\newcommand{\bvarepsilon}{\mbox{$\bf\varepsilon$}}

\newcommand{\dsize}{\displaystyle}
\renewcommand{\div}{\mathop{\rm div}\nolimits}
\newcommand{\ba}{\mbox{\boldmath{$a$}}}
\numberwithin{equation}{section}
\newcommand{\parent}[3]{\left #1 {#3} \right #2}
\newcommand{\graffe}[1]{\parent \{ \}{#1}}

\newcommand{\bignorm}[1]{\Bigl\lVert{#1}\Bigr\rVert} 
\newcommand{\norm}[1]{\lVert{#1}\rVert} 
\newcommand{\bfig}[2]{\begin{figure}\begin{center}\begin{picture}(341.8,#2)(
#1,0)}
\newcommand{\efig}[2]{\end{picture}\caption{#2.}\lbl{#1}\end{center}
\end{figure}}

\newcommand{\la}{\langle}
\newcommand{\ra}{\rangle}
\newcommand{\id}{{\boldsymbol{\mathbbm{1}}}}

\newtheorem{remark}{Remark}[section]

\newcommand{\be}{\begin{equation}}

\newcommand{\ee}{\end{equation}}


 \makeatletter
 \let\@fnsymbol\@arabic
 \makeatother

\begin{document}
\vskip-3truecm 
\title{
\vspace{-1.in} {\Large Existence result  for a dislocation based model of single crystal gradient plasticity with isotropic or linear kinematic hardening
\texttt{}}}

\author{{\large Fran\c{c}ois Ebobisse{\footnote{Corresponding author, Fran\c{c}ois Ebobisse,
Department of Mathematics and Applied Mathematics, University of
Cape Town, Rondebosch 7700, South Africa, e-mail:
francois.ebobissebille@uct.ac.za}} \quad and \quad  Patrizio Neff{\footnote
{Patrizio Neff, Lehrstuhl f\"ur  Nichtlineare Analysis und
Modellierung, Fakult\"{a}t f\"ur Mathematik, Universit\"at
Duisburg-Essen, Thea-Leymann Str. 9, 45127 Essen, Germany, e-mail:
patrizio.neff@uni-due.de, http://www.uni-due.de/mathematik/ag\b{
}neff}} \quad and\quad    Elias C. Aifantis{\footnote{Laboratory of Mechanics and Materials, Polytechnic School, Aristotle University of Thessaloniki, Thessaloniki GR-54124, Greece, e-mail: mom@mom.gen.auth.gr}} }\vspace{1mm}}



\date{\today}

 \maketitle

\begin{center}

\vspace{-1.5cm}\noindent

{\em In memory of Christian Miehe, 4.1.1956-14.8.2016, \\
 master of computational plasticity}\vspace{9mm}
 {{\small
\begin{abstract} We consider a dislocation-based rate-independent model of single crystal gradient plasticity with isotropic or linear kinematic hardening.  The model is weakly formulated through the so-called {\it primal form} of the flow rule as a variational inequality for which a result of existence and uniqueness is obtained using the functional analytical framework developed by Han-Reddy.

   \end{abstract}}}
\end{center}
\noindent {\bf Key words:} plasticity, single-crystal, gradient plasticity, variational modeling, dissipation function, defect energy, dislocation,  hardening, variational inequality.
 
\vskip.2truecm\noindent {\bf AMS 2010 subject classification:}
35D30, 35D35, 74C05, 74C15, 74D10, 35J25.
\newpage
 {\small \tableofcontents}

\section{Introduction}\label{Intro}

The goal of improving classical theories of plasticity which are incapable of modelling properly  material behaviour at small scales (meso/micron) as shown from experiments (\cite{FMAH1993,Fleckhutch1997,Stolkevans1998}), has led to the development of new theories which involve material length scales.  Since the pioneering work of Aifantis \cite{AIF1984}, where the yield-stress is
set to depend also on some derivative of a scalar measure of the
accumulated plastic distortion, there is now an abundant literature on models of gradient plasticity in the infinitesimal as well as in the finite deformation settings  for both polycrystaline and single crystal cases (see  \cite{AIF1984, AIF1987, MUHAIF1991,
Fleckhutch2001, AIF2003, GURT2002, GURT2004, GUD2004, GURTAN2005, GURTAN2005FD, GURTAN2009,
FleckWillisI2009, FleckWillisII2009, REDDY2011}). On the other hand, effort
has also been made in the past years to provide mathematical results
for the initial boundary values problems and inequalities describing
some models of (polycrystalline) gradient plasticity  (see for instance, \cite{DEMR1,
REM, EMR2008, NCA, EBONEFF, giacoluss,NESNEFF2012, NESNEFF2013,
ENR2015, EHN2016}). Several contributions on the computational aspects have
been made as well (\cite{DEMR2,NSW2009, BAREKLU2014, RWW2014}). However, mathematical results for models of single crystal gradient plasticity are rather scarce so far.

This paper deals with an infinitesimal strain single crystal plasticity formulation which accounts in a clear way for an energetic length scale effect. This is achieved by adding a quadratic measure of the dislocation density into the energy of cold work.

It is well accepted that moving dislocations are the carriers of plastic deformation. These moving dislocations leave the atomic lattice intact. In the idealized framework of a pure single crystal, the possible planes of these plastic deformations are assumed to be known in advance. The kinematics of plastic flow is thus constrained to these so-called {\it slip-planes}. The non-symmetric plastic distortion $p$ can therefore completely be described by the relation
\begin{equation}\label{dist-slip1}
p=\sum_{\alpha =1}^{n_{\mbox{\tiny slip}}}\gamma^\alpha\, l^\alpha\otimes\nu^\alpha\end{equation}
   where $l^\alpha$ and $\nu^\alpha$ are the $\alpha$-th slip direction  and the normal vector to the $\alpha$-th slip plane respectively, and $\gamma^\alpha$ describes the $\alpha$-th individual slip, for $\alpha=1,\ldots,n_{\mbox{\tiny slip}}$. The evolution law for the plastic distortion $p$ is then usually transfered to a system of evolution equations for the individual slips
   \begin{equation}\label{evol-slip1}
   \dot{\vectgam}=f(\Sigma_E(\vectgam))=\widehat{f}(\sigma,\vectgam)\,,\end{equation}
   where $\Sigma_E(\vectgam)$ is the  Eshelby-type stress tensor driving the evolution, $\sigma$ is the Cauchy stress tensor  and $f:\mathbb{R}^{3\times3}\to\mathbb{R}^{n_{\mbox{\tiny slip}}}$,  $\widehat{f}:\mathbb{R}^{3\times 3}\times\mathbb{R}^{n_{\mbox{\tiny slip}}}\to\mathbb{R}^{n_{\mbox{\tiny slip}}}$ are given functions driving the evolution. Here, we consider that not only the slips $\gamma^\alpha$ are intervening in the evolution for $\dot{\gamma}^\alpha$, but also gradients of slip $\nabla\dot{\gamma}^\alpha$. The rational for doing so is to incorporate a defect measure, i.e., a measure for the absence of compatibility of the plastic distortion. This  incompatibility is governed by the dislocation density tensor $\Curl p$ which, using equation (\ref{dist-slip1}) can be written (see \cite[Equ. (3.83)]{Han-ReddyBook})  as 
   \begin{equation}\label{curl-dist-slip}
   \Curl p=\sum_{\alpha =1}^{n_{\mbox{\tiny slip}}}(\nabla\gamma^\alpha\times\nu^\alpha)\otimes l^\alpha\,.\end{equation} 
   The idea in this paper is to specify a quadratic expression in $\Curl p$ such as 
     \begin{equation}\label{curl-dist-slip2}
   \frac12\,\mu\,L_c^2\,\norm{\Curl p}^2\,\end{equation} to augment the energy of cold work and to subsequently derive nonlocal evolution equations of the type   
   \begin{equation}\label{evol-slip2}
   \dot{\vectgam}=f(\sigma,\vectgam,\Curl\Curl p)=\widetilde{f}(\sigma,\vectgam,\nabla\vectgam,\nabla^2\vectgam)\end{equation}
   according to the guiding rule of maximal dissipation.
  
   In addition, we will consider local isotropic hardening as well as local kinematic hardening\footnote{We refer to ``local'' whenever these terms are already present in the classical theory without higher gradients.}. Moreover, in a special case we will try to bridge the gap to isotropic plasticity formulations by presenting the linear kinematic hardening through adding to the energy of cold work a term
    \begin{equation}\label{lin-kin-term1}
   \frac12\,\mu\,k_1\,\norm{\bvarepsilon_p}^2=\frac12\,\mu\,k_1\norm{\sym p}^2=\frac12\,\mu\,k_1\,\bignorm{\sum_\alpha\gamma^\alpha\sym( l^\alpha\otimes\nu^\alpha)}^2\,,
   \end{equation} reminiscent of Prager-type hardening.
 It turns out that while the combination of (\ref{curl-dist-slip2}) with (\ref{lin-kin-term1}) in an isotropic strain gradient context leads to a well-posed problem \cite{ENR2015}, in the single crystal setting (\ref{lin-kin-term1}) is, in general, not providing enough hardening for a mathematical treatment, even in the simplified case where the slip planes are mutually orthogonal. Only if we assume that the slip system are mutually orthogonal we can provide existence and uniqueness results in the linear kinematic setting as shown in Remark \ref{rem-choice-psi-kin}. The issue with   (\ref{lin-kin-term1}) is that it does  allow for free (infinitesimal) plastic rotations. Using instead of  (\ref{lin-kin-term1}) the term
 \begin{equation}\label{nosym-lin-kin}
  \frac12\,\mu\,k_1\norm{p}^2=\frac12\,\mu\,k_1\,\bignorm{\sum_\alpha\gamma^\alpha\, l^\alpha\otimes\nu^\alpha}^2\,,
   \end{equation} which is not invarient w.r.t. infinitesimal plastic rotations, would not cause any mathematical problem if the slip planes are assumed mutually orthogonal.
   
   In this sense, although single crystal plasticity needs the bureaucracy to treat simultaneously all slip systems $\mathbb{R}^{n_{\mbox{\tiny slip}}}$, the model is mathematically simpler than isotropic strain gradient plasticity models with plastic spin. These polycrystalline models need to be considered under additional invariance conditions which provide in effect less control for the plastic distortion. The long term goal remains to better understand the bridging between the computation of a single crystal and a polycrystalline sample.
   
   In this paper we do not consider ad hoc assumptions on the defect energy which would control the full slip gradient $\nabla\vectgam$, as this is not warranted by the theory. In our convex analytical setting we are, for the moment, bound to consider quadratic energies. We are aware of the fact that energies with linear growth like $\mu\,L_c\,\norm{\Curl p}^1$ represent important improvements of the modelling (see \cite{OHNOKUM2007}). This question will be postponed for further studies. In the next sections, we try to provide a descent, self-explaining introduction to the modelling aspects. Our focus is on transparency and clarity.

In \cite{GURT2002} a gradient theory of single-crystal plasticity that accounts for
the Burgers tensor and which is characterized for each slip system by a system of microforces which consists of a vector stress and a scalar force which are power conjugate to the corresponding slip in that system and its gradient respectively, is considered. The balance equations and boundary conditions in the model are derived through the virtual-power principle and the so-called {\it hard-slip conditions} are proposed. As shown in \cite{GURTNEED2005}, a discussion on suitable boundary conditions for the model proposed in \cite{GURT2002} depends on whether the Burger tensor, which is a fundamental ingredient in that theory, could be introduced directly in the virtual-power framework.
 Another type of boundary condition called {\it microhard condition} was then proposed in \cite{GURTNEED2005}.

The model in \cite{GURT2002, GURTNEED2005} with an isotropic hardening energy added to the free-nergy, was weakly formulated as a variational inequality in \cite{REDDY2011-II}
using the convex analytical framework developed for classical plasticity (polycrystalline and single crystal) in \cite{Han-ReddyBook}. Some well-posedness results were obtained as well in \cite{REDDY2011-II} depending on the choice of the defect energy density. Some computational results for this model are also obtained in \cite{ RWW2014}.
   
   The model in this paper presents some similarities and also some differences with the model in \cite{GURT2002, GURTNEED2005, REDDY2011-II}. In those papers, the energy of cold work is also augmented by   a quadratic expression with connections to the dislocation tensor $\Curl p$, which is rather expressed in each slip system in terms of 
   two other tensors representing the distributions of edge and screw dislocations on the given slip system. The defect energy is then  given as uncoupled and quadratic in those two new tensors. The evolution in  \cite[Equ. (2.28)]{REDDY2011-II} is governed in each slip system by both $\dot{\gamma}^\alpha$ and its gradient $\nabla\dot{\gamma}^\alpha$,  which then requires a further necessary condition on the micro-stresses and microforces driving the evolution  in that slip system. In our model, the evolution  is only governed in each slip system  by $\dot{\gamma}^\alpha$ and we do not need to involve the distributions of edge and screw dislocations on each slip system for the mathematical well-posedness.
   
   This paper is mainly divided into two parts. We first derive in section \ref{single-iso} both strong and weak formulations of  the model with isotropic hadening then we study its well-posedness. Section \ref{single-kin}
 is devoted to the model with linear kinematic hardening.
 \section{Some notational agreements and
definitions}\label{Notations} Let $\Omega$ be a bounded domain
 in $\BBR^3$ with Lipschitz continuous boundary $\partial\Omega$, which is occupied by the elastoplastic
body in its undeformed configuration. Let $\Gamma_{\mbox{\scriptsize D}}$ be a smooth
subset of $\partial\Omega$ with non-vanishing $2$-dimensional
Hausdorff measure. A material point in $\Omega$ is denoted by $x$
and the time domain under consideration is the interval $[0,T]$.\\
 For every $a,\,b\in\BBR^3$, we let $\la a,b\ra_{\BBR^3}$ denote the scalar
 product on $\BBR^3$ with associated vector
norm $\norm{a}^2_{\BBR^3} = \la a, a\ra_{\BBR^3}$. We denote by
$\BBR^{3\times 3}$ the set of real $3\times 3$ tensors. The standard
Euclidean scalar product on $\BBR^{3\times 3}$ is given by $\la
A,\,B\ra_{\BBR^{3\times 3}} = \mbox{tr}\,\bigl[AB^T\bigr]$, where
$B^T$ denotes the transpose tensor of $B$. Thus, the Frobenius
tensor norm is $\norm{A}^2_{\mathbb{R}^{3\times 3}} = \la A,\,A\ra_{\BBR^{3\times 3}}$.
 In the following we omit the subscripts $\BBR^3$ and $\BBR^{3\times 3}$. The identity tensor on $\BBR^{3\times 3}$ will be denoted by
  $\id$, so that $\mbox{tr}(A) = \la A, \id\ra$.  The set
$\so(3):=\{X\in\BBR^{3\times 3}\,|\,\,X^T=-X\}$ is the Lie-Algebra
of skew-symmetric tensors.
 We let
$\mbox{Sym\,}(3):=\{X\in\BBR^{3\times 3}\,|\,\,X^T=X\}$ denote the
 vector space of symmetric tensors and $\sL(3):=\{X\in\BBR^{3\times
3}\,|\,\,\mbox{tr\,}(X)=0\}$ be the Lie-Algebra of traceless
tensors. For every $X\in\BBR^{3\times 3}$, we set
$\sym(X)=\frac12\bigl(X+X^T\bigr)$,
$\skew\,(X)=\frac12\bigl(X-X^T\bigr)$ and
$\dev(X)=X-\frac13\mbox{tr}\,(X)\,\id\in\sL(3)\,$ for the symmetric
part, the skew-symmetric part and the deviatoric part of $X$,
respectively. Quantities which are constant in space will be denoted
with an overbar, e.g., $\overline{A}\in\so(3)$ for the function
$A:\mathbb{R}^3\to\so(3)$ which is constant with constant value
$\overline{A}$.

The body is assumed to undergo infinitesimal deformations. Its
behaviour is governed by a set of equations and constitutive
relations. Below is a list of variables and parameters used
throughout the paper with their significations:\begin{itemize}
\item[$\bullet$] $u$  is the displacement of the macroscopic material
points;

\item[$\bullet$] $p$ is the infinitesimal plastic distortion variable which is a
non-symmetric second order tensor, incapable of sustaining
volumetric changes; that is, $p\in\sL(3)$. The tensor $p$\,
represents the average plastic slip; $p$ is not a state-variable, while the rate $\dot{p}$ is;

\item[$\bullet$] $e=\nabla u -p$ is  the infinitesimal elastic distortion which  is in general a
non-symmetric second order tensor and is a state-variable;

\item[$\bullet$] $\bvarepsilon_p=\sym p$ is the symmetric infinitesimal plastic strain
tensor, which is  trace free, {$\bvarepsilon_p\in\sL(3)$;} $\bvarepsilon_p$ is not a state-variable;
the rate $\dot{\bvarepsilon}_p$ is a state-variable;

\item[$\bullet$] $\bvarepsilon_e=\sym\nabla u -\bvarepsilon_p$ is the symmetric infinitesimal  elastic
strain tensor and is a state-variable;

\item[$\bullet$] $\sigma$   is the Cauchy stress tensor which is a symmetric
second order tensor and is a state-variable;

\item[$\bullet$] $\yieldzero$ is the initial
yield stress for plastic variables $p$ or $\bvarepsilon_p:=\sym p$ and is a state-variable;

\item[$\bullet$] $f$ is the body force;

\item[$\bullet$] $\Curl p=\alpha$ is the dislocation density
tensor satisfying the so-called Bianchi identities $\Div\alpha=0$ and is a state-variable;


\item[$\bullet$] $\gamma^\alpha$ is the slip in the $\alpha$-th slip system in single crystal plasticity while $l^\alpha$ is the slip direction and $\nu^\alpha$ is the normal vector to the slip plane with $\alpha=1,\ldots,n_{\mbox{\tiny slip}}$. Hence, $p=\dsize\sum_\alpha\gamma^\alpha\,l^\alpha\otimes\nu^\alpha$.
\end{itemize}
\vskip.2truecm\noindent For isotropic media, the fourth order
isotropic elasticity tensor $\C_{\mbox{\scriptsize
iso}}:\mbox{Sym}(3)\to\mbox{Sym}(3)$ is given by
\begin{equation}
\C_{\mbox{\scriptsize iso}}\sym X = 2\mu\,\dev\,\sym X+\kappa
\,\tr(X) \id =2\mu\,\sym X+\lambda\,\tr(X)\id\label{C}
\end{equation}
for any second-order tensor $X$, where $\mu$ and $\lambda$ are the
Lam{\'e} moduli satisfying
\begin{equation}\label{Lame-moduli}
\mu>0\quad\mbox{ and }\quad 3\lambda +2\mu>0\,,
\end{equation} and $\kappa>0$ is the bulk modulus.
These conditions suffice for pointwise positive definiteness of the
elasticity tensor in the sense that there exists a constant $m_0 >
0$ such that
\begin{equation}
\forall X\in\mathbb{R}^{3\times 3}\mbox{ :}\,\quad\la \sym X,\C_{\mbox{\scriptsize iso}}\sym X\ra \geq m_0\,
\lVert\sym X\rVert^2\,. \label{ellipticityC}
\end{equation}

\vskip.2truecm\noindent The space of square integrable functions is
$L^2(\Omega)$, while the Sobolev spaces used in this paper are:
\begin{eqnarray}\label{sobolev-spaces}
\nonumber \mbox{H}^1(\Omega)&=&\{u\in
L^2(\Omega)\,\,|\,\,\mbox{grad\,}u\in
L^2(\Omega)\}\,,\qquad\quad\mbox{ grad}\,=\nabla\,,\\
\nonumber
&&\norm{u}^2_{H^1(\Omega)}=\norm{u}^2_{L^2(\Omega)}+\norm{\mbox{grad\,}u}^2_{L^2(\Omega)}\,,\qquad\forall u\in\mbox{H}^1(\Omega)\,,\\
 \mbox{H}(\mbox{curl};\Omega)&=&\{v\in
L^2(\Omega)\,\,|\,\,\mbox{curl\,}v\in
L^2(\Omega)\}\,,\qquad\mbox{curl\,}=\nabla\times\,,\\
\nonumber &&\norm{v}^2_{\mbox{\scriptsize
H}(\mbox{curl};\Omega)}=\norm{v}^2_{L^2(\Omega)}+\norm{\mbox{curl\,}v}^2_{L^2(\Omega)}\,,\,\,\quad\forall
v\in\mbox{H}(\mbox{curl;\,}\Omega)\,.
\end{eqnarray}
For every $X\in C^1(\Omega,\,\BBR^{3\times 3})$ with rows
$X^1,\,X^2,\,X^3$, we use in this paper the definition of $\Curl X$
in \cite{NCA, SVEN}:
\begin{equation}\label{def-Curl}\Curl X =\left(\begin{array}{l}\mbox{curl\,}X^1\,\,-\,\,-\\
\mbox{curl\,}X^2\,\,-\,\,-\\
\mbox{curl\,}X^3\,\,-\,\,-\end{array}\right)\in\BBR^{3\times
3}\,,\end{equation} for which $\Curl\,\nabla v=0$ for every $v\in
C^2(\Omega,\,\BBR^3)$. Notice that the definition of $\Curl X$
above is such that $(\Curl X)^Ta=\mbox{curl\,}(X^Ta)$ for every
$a\in\BBR^3$ and this clearly corresponds to the transpose of the
Curl of a tensor as defined in
\cite{GURTAN2005, GURTAN-BOOK}.\\

The following function spaces and norms will also be used later.
\begin{eqnarray}\label{Curl-spaces}
\nonumber \mbox{H}(\mbox{Curl};\,\Omega,\,\BBR^{3\times
3})&=&\Bigl\{X\in L^2(\Omega,\,\BBR^{3\times
3})\,\,\bigl|\,\,\mbox{Curl\,}X\in
L^2(\Omega,\,\BBR^{3\times 3})\Bigr\}\,,\\
\norm{X}^2_{\mbox{\scriptsize H}(\mbox{\scriptsize
Curl};\Omega)}&=&\norm{X}^2_{L^2(\Omega)}+\norm{\mbox{Curl\,}X}^2_{L^2(\Omega)}\,,\quad\forall
X\in\mbox{H}(\mbox{Curl;\,}\Omega,\,\BBR^{3\times 3})\,,\\
\nonumber
\mbox{H}(\mbox{Curl};\,\Omega,\,\mathbb{E})&=&\Bigl\{X:\Omega\to\mathbb{E}\,\,\bigl|\,\,X\in
\mbox{H}(\mbox{Curl};\,\Omega,\,\BBR^{3\times 3})\Bigr\}\,,
\end{eqnarray}
for $\mathbb{E}:=\sL(3)$ or
$\mbox{Sym}\,(3)\cap\sL(3)$.\vskip.2truecm\noindent
 We also consider the space
\begin{equation}\label{spacep-bc}\mbox{H}_0(\mbox{Curl};\,\Omega,\,\Gamma_{\mbox{\scriptsize D}},\BBR^{3\times 3})\end{equation} as
the completion in the norm in  (\ref{Curl-spaces}) of the space
$\bigl\{X\in C^\infty(\Omega,\,\BBR^{3\times
3})\,\,\bigl|\,\,\,X\times\,n|_{\Gamma_{\mbox{\scriptsize D}}}=0\bigr\}\,.$ Therefore, this
space generalizes the tangential Dirichlet boundary condition
$$X\times\,n|_{\Gamma_{\mbox{\scriptsize D}}}=0\,$$
to be satisfied by the plastic distortion $p$. The space
$$\mbox{H}_0(\mbox{Curl};\,\Omega,\,\Gamma_{\mbox{\scriptsize D}},\mathbb{E})$$ is defined
as
 in (\ref{Curl-spaces}). \vskip.2truecm\noindent
 The divergence operator Div on second order
tensor-valued functions is also defined row-wise as
\begin{equation}\label{def-div}\mbox{Div}\,X=\left(\begin{array}{l}\mbox{div\,}X_1\\
\mbox{div\,}X_2\\
\mbox{div\,}X_3\end{array}\right)\,.\end{equation}

\section{The kinematics  of single crystal plasticity}
Single-crystal plasticity is based on the assumption that the plastic deformation happens through crystallograpic shearing which represents the dislocation motion along  specific slip systems,  each being characterized
by a plane with unit normal $\nu^\alpha$ and slip direction $l^\alpha$ on that plane, and slips $\gamma^\alpha$ ($\alpha = 1, \ldots ,n_{\mbox{\tiny slip}}$). The flow rule
for the plastic distortion $p$ is written at the slip system level by means of the orientation tensor $m^\alpha$ defined as
\begin{equation}\label{orient-tensor} m^\alpha:=l^\alpha\otimes \nu^\alpha\,.\end{equation}
Under these conditions the plastic distortion $p$  takes the form
\begin{equation}\label{orientens} 
p=\sum_{\alpha=1}^{n_{\mbox{\tiny slip}}}\gamma^\alpha \,m^\alpha\end{equation}
so that the plastic strain $\bvarepsilon_p=\sym p$ is
\begin{equation}\label{orientens2}\bvarepsilon_p=\sum_{\alpha=1}^{n_{\mbox{\tiny slip}}} \gamma^\alpha \sym(m^\alpha)=\frac12\sum_{\alpha=1}^{n_{\mbox{\tiny slip}}}      \gamma^\alpha(l^\alpha\otimes \nu^\alpha+\nu^\alpha\otimes l^\alpha)\end{equation} and $\tr(p)=\tr(\bvarepsilon_p)=0$ since by assumption $l^\alpha\perp \nu^\alpha$.\vskip.2truecm\noindent
For the slips $\gamma^\alpha$ ($\alpha=1,\dots,n_{\mbox{\tiny slip}}$) we set
$$\mbox{\underline{$\gamma$}}:=(\gamma^1,\ldots,\gamma^{n_{\mbox{\tiny slip}}})\,.$$
Therefore, we get from (\ref{orientens2}) that
\begin{equation}\label{3rd-order-m}
p=\overline{m}\,\vectgam\,,\end{equation} where $\overline{m}$ is the third order tensor defined as \begin{equation}\label{orientens3} \overline{m}_{ij\alpha}:=m^\alpha_{ij}=l^\alpha_i\nu^\alpha_j\quad\mbox{ for}\quad i,j=1,2,3\, \mbox{ and } \alpha=1,\ldots,n_{\mbox{\tiny slip}}\,.\end{equation}  We will extensively make use of the identity (\ref{3rd-order-m}). \vskip.2truecm\noindent
Let $\mbox{\underline{$\eta$}}:=(\eta^1,\ldots,\eta^{n_{\mbox{\tiny slip}}})$ with $\eta^\alpha$ being a hardening variable in the $\alpha$-th slip system.

\section{The description of the model with isotropic hardening}\label{single-iso}
\subsection{The balance equation} The conventional macroscopic
force balance leads to the equation of equilibrium
\begin{equation}
\Div \sigma + f = \bzero \label{equil}
\end{equation}
in which $\sigma$ is the infinitesimal symmetric Cauchy stress and
$f$ is the body force. 
\subsection{Constitutive equations.} The
constitutive equations are obtained from a free energy imbalance
together with a flow law that characterizes plastic behaviour. Since
the model under study involves plastic spin by which we mean that
the plastic distortion $p$ is not symmetric, we consider directly an additive
decomposition of the displacement gradient $\nabla u$ into elastic
and plastic components $e$ and $p$, so that
\begin{equation}
\nabla u = e + p\,, \label{displ-grad}
\end{equation}
with the nonsymmetric plastic distortion $p$ incapable of
sustaining volumetric changes; that is,
\begin{equation}
\tr(p)=\tr(\sym p)=\tr(\bvarepsilon_p) = 0\,. \label{trEp}
\end{equation}
 Here, $\bvarepsilon_e=\sym e=\sym (\nabla u-p)$ is the
infinitesimal elastic strain and $\bvarepsilon_p=\sym p$ is the
plastic strain while $\sym \nabla u=(\nabla
u+\nabla u^T)/2$ is the total strain.
 \subsubsection{The free-energy}\label{free-eng-iso}
  We consider  a free
energy in the additively separated form
\begin{eqnarray}\label{free-eng}
\Psi(\nabla u,\vectgam,\Curl p,\vecteta):
&=&\underbrace{\Psi^{\mbox{\scriptsize lin}}_e(\sym
e)}_{\mbox{\small elastic energy}}\,\,
+\,\,\,\underbrace{\Psi^{\mbox{\scriptsize lin}}_{\mbox{\scriptsize
curl}}(\Curl p)}_{\mbox{\small defect
energy (GND)}}\\
\nonumber
&&\qquad\quad+\hskip-1.truecm\,\,\underbrace{\,\,\Psi_{\mbox{\scriptsize
 iso}}(\vecteta)}_{\begin{array}{c}\mbox{\small isotropic}\\
\mbox{\small hardening energy (SSD)}\end{array}}\,,
 \end{eqnarray} where
 \begin{equation}\label{free-eng-expr}\left\{\begin{array}{lcl} \Psi^{\mbox{\scriptsize
lin}}_e(\sym e): &=&\frac12\,\la\sym e,\C_{\mbox{\scriptsize iso}} \sym
e\ra,\\\\
 \Psi^{\mbox{\scriptsize lin}}_{\mbox{\scriptsize
curl}}(\Curl p):&=&\frac12\,\mu\, L_c^2\,\norm{\Curl p}^2=\frac12\,\mu\,L_c^2\norm{\Curl(\overline{m}\,\underline{\gamma})}^2\,,\\
\\ \Psi_{\mbox{\scriptsize
 iso}}(\vecteta) &:=&\frac12\,\mu\,k_2\,\norm{\vecteta}^2\,=\,\frac12\mu\,k_2\,\sum_\alpha|\eta^\alpha|^2
 \,.\end{array}\right.\end{equation}
Here,  $L_c\geq0$ is an energetic length scale which characterizes the contribution of the defect energy density to the system,  $k_2$
a positive nondimensional isotropic hardening constant. The
defect energy is conceptually related to geometrically necessary
dislocations (GND). It is formed by the long-ranging stress-fields of excess dislocations and may be recovered by appropriate inelastic deformation.
The isotropic hardening energy is related to
statistically stored dislocations (SSD).

 \subsubsection{The
derivation of the dissipation inequality}\label{dis-ineq-iso}
The local free-energy imbalance states that
\begin{equation}
\dot{\Psi} - \la\sigma,\dot{e}\ra - \la\sigma,\dot{p}\ra  \leq 0\
. \label{2ndlaw}
\end{equation}
Now we expand the first term, substitute (\ref{free-eng}) and get
\begin{equation}\label{exp-2ndlaw-crystal1}
\la\C_{\mbox{\scriptsize iso}}\sym
e-\sigma,\sym\dot{e}\ra-\la\dev\sigma,\overline{m}\,\underline{\dot{\gamma}}\ra+\mu\,k_2\,\la\vecteta,\dot{\vecteta}\ra+\mu
L_c^2\la\Curl(\overline{m}\,\underline{\gamma}),\Curl(\overline{m}\,\underline{\dot{\gamma}})\ra\leq0\,,
\end{equation}
Since the inequality (\ref{exp-2ndlaw-crystal1}) must be satisfied for whatever elastic-plastic deformation mechanism, inlcuding purely elastic ones (for which  $\underline{\dot{\gamma}}=\dot{\vecteta}=0$), inequality (\ref{exp-2ndlaw-crystal1}) implies the usual infinitesimal elastic stress-strain relation
\begin{eqnarray}
 \sigma \,=\, \C_{\mbox{\scriptsize iso}}\,\varepsilon_e=2\mu\, \sym(\nabla
u-p)+\lambda\, \tr(\nabla u-p)\id \,=\,2\mu\, (\sym(\nabla
u)-\bvarepsilon_p)+\lambda\, \tr(\nabla u)\id\label{elasticlaw}
\end{eqnarray}
and the local reduced dissipation inequality
\begin{equation}\label{dissip-ineq-crystal1}
\la\dev\sigma,\overline{m}\,\underline{\dot{\gamma}}\ra-\mu\,k_2\,\la\vecteta,\dot{\vecteta}\ra-\mu
L_c^2\la\Curl(\overline{m}\,\underline{\gamma}),\Curl(\overline{m}\,\underline{\dot{\gamma}})\ra\geq0\,,
\end{equation} that we integrate over $\Omega$ and get
\begin{eqnarray}
\nonumber\hskip-4truecm0&\leq&\int_{\Omega}\Bigl[\la\dev\sigma,\overline{m}\,\underline{\dot{\gamma}}\ra-\mu\,k_2\,\la\vecteta,\dot{\vecteta}\ra-\mu
L_c^2\la\Curl(\overline{m}\,\underline{\gamma}),\Curl(\overline{m}\,\underline{\dot{\gamma}})\ra\Bigr]\,dx\\
 \nonumber\hskip-4truecm &=&\int_{\Omega}\Bigl[\la\dev\sigma,\overline{m}\,\underline{\dot{\gamma}}\ra-\mu\,k_2\,\la\vecteta,\dot{\vecteta}\ra-\mu
L_c^2\la\Curl\Curl(\overline{m}\,\underline{\gamma}),\overline{m}\,\underline{\dot{\gamma}}\ra\Bigr]\,dx\\
 \nonumber &&\hskip4truecm-\,\,\sum_{i=1}^3\mu\,L_c^2\int_\Omega\div\left((\overline{m}\,\underline{\dot{\gamma}})^i\times(\Curl(\overline{m}\,\underline{\gamma}))^i\right)\,dx\\
\hskip-4truecm &=&\int_{\Omega}\Bigl[\la\dev\sigma,\overline{m}\,\underline{\dot{\gamma}}\ra-\mu\,k_2\la\vecteta,\dot{\vecteta}\ra-\mu
L_c^2\la\Curl\Curl(\overline{m}\,\underline{\gamma}),\overline{m}\,\underline{\dot{\gamma}}\ra\Bigr]\,dx\\
 \nonumber &&\hskip4truecm+\,\,\sum_{i=1}^3\mu\,L_c^2\int_{\partial\Omega}\la(\overline{m}\,\underline{\dot{\gamma}})^i\times(\Curl(\overline{m}\,\underline{\gamma}))^i,n\ra\,dS\,.
\end{eqnarray}
In order to obtain a dissipation
inequality in the spirit of classical plasticity, we assume that the
infinitesimal plastic distortion $p$ (and hence the slips through $\overline{m}\,\underline{\gamma}$) satisfies the so-called {\it
linearized insulation condition}
\begin{equation}\label{lin-sul-iso}\sum_{i=1}^3\mu\,L_c^2\int_{\partial\Omega}\la(\overline{m}\,\underline{\dot{\gamma}})^i\times(\Curl(\overline{m}\,\underline{\gamma}))^i,n\ra\,dS=0\,.\end{equation} Under (\ref{lin-sul})  we then
obtain a global version of the reduced dissipation inequality
\begin{equation}\label{global-dissip-ineq1}\int_{\Omega}\Bigl[\la\dev\sigma,\overline{m}\,\underline{\dot{\gamma}}\ra-\mu\,k_2\,\la\vecteta,\dot{\vecteta}\ra-\mu
L_c^2\la\dev\Curl\Curl(\overline{m}\,\underline{\gamma}),\overline{m}\,\underline{\dot{\gamma}}\ra\Bigr]\,dx\geq0\,.\end{equation}
That is,
\begin{eqnarray}\label{global-dissip-ineq2}
\int_\Omega\sum_\alpha\Bigl[\la\dev\sigma,m^\alpha\ra\,\dot{\gamma}^\alpha\ra-\mu\,k_2\,\eta^\alpha\,\dot{\eta}^\alpha-\la\mu\,L_c^2\la\dev\Curl\Curl(\overline{m}\,\underline{\gamma}),m^\alpha\ra\,\dot{\gamma}^\alpha\Bigr]\,dx\geq0\,.\end{eqnarray} Hence, we get
 \begin{equation}\label{dissip-ineq-final}
 \int_\Omega\sum_\alpha\Bigl[\tau^\alpha_{\mbox{\tiny E}}\,\dot{\gamma}^\alpha+g^\alpha\,\dot{\eta}^\alpha\Bigr]\,dx\geq0\,\end{equation}
 where we set
\begin{eqnarray}\tau^\alpha_{\mbox{\tiny E}}&:=&\tau^\alpha +s^\alpha_{\mbox{\tiny nonloc}}\\
\mbox{ }\hskip-1truecm \tau^\alpha&:=&\la\dev\sigma,m^\alpha\ra\quad\mbox{\,\,\,\,(resolved shear stress for the $\alpha$-th slip system)}\,,\\
s^\alpha_{\mbox{\tiny nonloc}}&:=&-\,\mu\,L_c^2\,\la\dev\Curl\Curl(\overline{m}\,\underline{\gamma}),m^\alpha\ra\,,\\
\nonumber&&\hskip1.8truecm\mbox{\,\,(nonlocal backstress for the $\alpha$-th slip system)},\\
g^\alpha:&=&-\,\mu\,k_2\,\eta^\alpha\quad\mbox{(thermodynamic force power-conjugate to $\dot{\eta}^\alpha$)} \,.
\end{eqnarray} 
 Notice that 
 $$\tau^\alpha_{\mbox{\tiny E}}=\la\dev\Sigma_{\mbox{\tiny E}},m^\alpha\ra$$ 
 with  $\Sigma_E$  being the non-symmetric Eshelby-type stress tensor defined by
\begin{equation}\label{Eshelby-crystal} \Sigma_E:=\sigma-\mu\,L_c^2\,\Curl\Curl(\overline{m}\,\underline{\gamma})\,.
\end{equation}
 The local reduced dissipation inequality can also be written in compact form as 
\begin{equation}\label{reduced-diss-comp}
\sum_\alpha\la\Sigma_p^\alpha,\dot{\Gamma}_p^\alpha\ra\geq0\,,
\end{equation}where we define 
\begin{equation}\label{gen-var}\Sigma_p^\alpha:=(\tau^\alpha_E,g^\alpha)\quad\mbox{ and }\quad \Gamma_p^\alpha:=(\gamma^\alpha,\eta^\alpha)\,.
\end{equation}
\subsubsection{The boundary conditions on the plastic
distortion}\label{bc-distortion}
 The condition (\ref{lin-sul}) is satisfied if
we
 assume for instance that the boundary is a perfect conductor. This means that the tangential component of $p=\overline{m}\,\underline{\gamma}$ vanishes on
$\partial\Omega$.
 In the context of dislocation dynamics these conditions express the requirement
  that there is no flux of the Burgers vector across a hard boundary.
Gurtin \cite{GURT2004} and also Gurtin and Needleman \cite{GURTNEED2005} introduce the
following different types of boundary conditions for the plastic distortion%
\begin{align}
    \dot{p}\times\,n|_{\Gamma_{\rm hard}}&=0 \quad \text{"micro-hard" (perfect conductor)} \notag \\
    \dot{p}|_{\Gamma_{\rm hard}}&=0 \quad \text{"hard-slip"}\quad\mbox{(in the context of crytal plasticity)}\\
    (\Curl \dot{p})\times\,n|_{\Gamma_{\rm hard}}&=0 \quad \text{"micro-free"}\, . \notag
\end{align}
Notice that ``hard-slip'' boundary condition $\dot{p}|_{\Gamma_{\rm hard}}=0$ makes more sense in the context where the defect energy is formulated in terms of $\nabla p$. In our current situation
 we specify a sufficient condition for the micro-hard boundary
condition, namely \begin{equation}\label{bc-plastic}
       p\times\,n|_{\Gamma_{\rm hard}}=(\overline{m}\,\underline{\gamma})\times n|_{\Gamma_{\rm hard}}=0
\end{equation}
and assume for simplicity that $\Gamma_{\rm
hard}=\Gamma_{\mbox{\tiny D}}$. Note that this boundary condition
constrains the plastic slip in tangential direction only, which is
what we expect to happen at the physical boundary $\Gamma_{\rm
hard}$.
\subsubsection{The flow rule}\label{flow-law}
We consider a yield function on the $\alpha$-th slip system
defined by
\begin{equation}
\label{yield-funct-kinematic} \phi(\Sigma^\alpha_p):= |\tau^\alpha_{\mbox{\tiny E}}|+g^\alpha - \yieldzero\quad\mbox{ for }\quad\Sigma^\alpha_p=(\tau^\alpha_{\mbox{\tiny E}},g^\alpha)\,.
\end{equation}
Here,  $\yieldzero$ is the initial yield stress of the material, that we assume to be constant on all slip systems and therefore, $\yieldlimit^\alpha:=\sigma_0-g^\alpha$ represents the current yield stress for the $\alpha$-th slip system{\footnote{Note that, for the sake of simplicity, the presented isotropic hardening
rule $g^\alpha$ does not involve latent hardening and the associated interaction matrix, see \cite{franciosizaoui91}
for a discussion on uniqueness in the presence of latent hardening.}. So the set
of admissible  generalized stresses for the $\alpha$-th slip system is defined as
\begin{equation}\label{admiss-stress-kin}
\mathcal{K}^\alpha:=\left\{\Sigma^\alpha_p=(\tau^\alpha_{\mbox{\tiny
E}},g^\alpha)\,\,|\,\,\phi(\Sigma_p^\alpha)
   \leq0,\,\,g^\alpha\leq0\right\}\,,\end{equation} with its interior $\mbox{Int}(\mathcal{K}^\alpha)$ and its boundary $\partial\mathcal{K}^\alpha$ being the generalized elastic region and the yield surface for the $\alpha$-th slip system, respectively.
   \vskip.1truecm\noindent
The  principle of maximum dissipation associated with the $\alpha$-th slip system gives  us the normality
law\begin{equation}\label{normalcone} \dot{\Gamma}_p^\alpha\in
N_{\mathcal{K}^\alpha}(\Sigma^\alpha_p)\,,\end{equation}
where $\dsize N_{\mathcal{K}^\alpha}(\Sigma^\alpha_p)$
denotes the normal cone to $\mathcal{K}^\alpha$ at
$\dsize\Sigma_p^\alpha$. That is,  $\dot{\Gamma}^\alpha_p$ satisfies
\begin{equation}
\la\overline{\Sigma}^\alpha - \Sigma^\alpha_p,\dot{\Gamma}^\alpha_p\ra \leq
0\ \quad \mbox{for all}\ \overline{\Sigma}^\alpha \in {\mathcal{K}^\alpha}\ .
\label{normality2-alpha}
\end{equation} Notice that $N_{\mathcal{K}^\alpha}=\partial\Chi_{\mathcal{K}^\alpha}$, where
$\Chi_{\mathcal{K}^\alpha}$ denotes the indicator function of the set $\mathcal{K}^\alpha$
and $\partial\Chi_{\mathcal{K}^\alpha}$ denotes the subdifferential of the function
$\Chi_{\mathcal{K}^\alpha}$.\\ Whenever the yield surface $\partial\mathcal{K}^\alpha$ is
smooth at $\dsize\Sigma_p^\alpha$ then
$$\dot{\Gamma}^\alpha_p\in
N_{\mathcal{K}^\alpha}(\Sigma^\alpha_p)\quad\Rightarrow\quad\exists\lambda^\alpha\mbox{ such that }
\dot{\gamma}^\alpha=\lambda^\alpha\,\frac{\tau^\alpha_{\mbox{\tiny E}}}{|\tau^\alpha_{\mbox{\tiny E}}|}\quad\mbox{and}\quad\dot{\eta}^\alpha=\lambda^\alpha=|\dot{\gamma}^\alpha|$$ with the Karush-Kuhn
Tucker conditions: $\lambda^\alpha\geq0$, $\phi(\Sigma^\alpha_p)\leq0$ and $\lambda^\alpha\,\phi(\Sigma^\alpha_p)=0$\,.\\
Using convex analysis (Legendre-transformation) we find that
\begin{eqnarray}\label{dualflowlaw-crystal} &\underbrace{\dot{\Gamma}^\alpha_p\in
\partial\Chi_{\mathcal{K}^\alpha}(\Sigma^\alpha_p)}_{\mbox{\bf flow rule in its dual formulation for the $\alpha$-th slip system}}&\\
\nonumber&\Updownarrow& \\
&\underbrace{\Sigma^\alpha_p\in
\partial \Chi^*_{\mathcal{K}^\alpha} (\dot{\Gamma}^\alpha_p)\,}_{\mbox{\bf flow rule in its primal formulation  for the $\alpha$-th slip system}}&\label{primalflowlaw-crystal}
\end{eqnarray} where $\Chi^*_{\mathcal{K}^\alpha}$ is the Fenchel-Legendre dual of the function $\Chi_{\mathcal{K}^\alpha}$ denoted in this context by $\mathcal{D}^\alpha_{\mbox{\scriptsize iso}}$,
 the one-homogeneous dissipation function for the $\alpha$-th slip system. That is, for every $\Gamma^\alpha=(q^\alpha,\beta^\alpha)$,
\begin{eqnarray}\label{dissp-function-iso-slip}
\nonumber\mathcal{D}^\alpha_{\mbox{\scriptsize iso}}(\Gamma^\alpha)
&=&\sup\graffe{\la\Sigma^\alpha_p,\Gamma^\alpha\ra\,\,|\,\,\Sigma^\alpha_p\in\mathcal{K}^\alpha}\\
\nonumber&=&
\sup \graffe{\tau^\alpha_E\,q^\alpha +g^\alpha\beta^\alpha \,\,|\,\
\phi(\Sigma^\alpha_E,g^\alpha)\leq0,\,\,g^\alpha\leq0}\\&=&\left\{\begin{array}{ll}
      \yieldzero\,|q^\alpha| &\mbox{ if } |q^\alpha|\leq\beta^\alpha\,,\\
      \infty &\mbox{ otherwise.}\end{array}\right.
\end{eqnarray} We get from the definition of the subdifferential ($\Sigma^\alpha_p \in
\partial \Chi^*_{\mathcal{K}^\alpha} (\dot{\Gamma}^\alpha_p)$) that
\begin{equation}
\mathcal{D}^\alpha_{\mbox{\scriptsize iso}} (\Gamma^\alpha) \geq \mathcal{D}^\alpha_{\mbox{\scriptsize iso}}(\dot{\Gamma}^\alpha_p) +
\la\Sigma^\alpha_p,\Gamma^\alpha-\dot{\Gamma}^\alpha_p\ra\quad \mbox{for any
} \Gamma^\alpha.\label{dissineq-slip}\end{equation}
That is,
\begin{equation} \mathcal{D}^\alpha_{\mbox{\scriptsize iso}}(q^\alpha,\beta^\alpha)\geq \mathcal{D}^\alpha_{\mbox{\scriptsize iso}}(\dot{\gamma}^\alpha,\dot{\eta}^\alpha)+\tau^\alpha_{\mbox{\tiny E}}\,(q^\alpha-\dot{\gamma}^\alpha)+g^\alpha(\beta^\alpha-\dot{\eta}^\alpha)\quad \mbox{for any
} (q^\alpha,\beta^\alpha).\label{dissinequality-slip2}\end{equation}
\vskip.3truecm\noindent In the next sections, we present a complete mathematical analysis of the model including  both strong and weak formulations as well as a corresponding
existence result.

\subsection{Mathematical analysis of the model}\label{math-anal-iso}
\subsubsection{Strong formulation}\label{strong-iso-crystal} 
\label{strong-iso-crystal} To summarize, we have obtained the following strong
formulation for the model of single crystal infinitesimal gradient plasticity with
 isotropic  hardening.  Given $f\in \SFH^1(0,T;L^2(\Omega,\mathbb{R}^3))$, the goal is to find:
\begin{itemize}\item[(i)] the displacement $u\in \SFH^1(0,T; H^1_0(\Omega,{\Gamma_{\mbox{\scriptsize D}}},\mathbb{R}^3))$,
\item[(ii)] the infinitesimal plastic slips $\gamma^\alpha\in
\SFH^1(0,T;L^2(\Omega))$ for $\alpha=1,\ldots,n_{\mbox{\tiny slip}}$  with\\ $\Curl(\overline{m}\,\underline{\gamma})\in\SFH^1(0,T;H(\Curl;\Omega, \mathbb{R}^{3\times 3})$, 
\item[(iii)] the hardening variables $\eta^\alpha\in
\SFH^1(0,T;L^2(\Omega))$ for $\alpha=1,\ldots,n_{\mbox{\tiny slip}}$, \end{itemize}
 such that the content of Table \ref{table:isohard-crystal} holds.
 
 \begin{table}[ht!]{\footnotesize\begin{center}
\begin{tabular}{|ll|}\hline &\qquad\\
 {\em Additive split of distortion:}& $\nabla u =e +p$,\quad $\bvarepsilon_e=\mbox{sym}\,e$,\quad $\bvarepsilon_p=\sym p$\\
{\em Plastic distortion in slip system:}&$p=\dsize\sum_{\alpha=1}^{n_{\mbox{\tiny silp}}}      \gamma^\alpha\,m^\alpha$ with $m^\alpha= l^\alpha\otimes \nu^\alpha$, \quad $\tr(p)=0$\\
{\em Equilibrium:} & $\mbox{Div}\,\sigma +f=0$ with
$\sigma=\C_{\mbox{\tiny iso}}\bvarepsilon_e=\C_{\mbox{\tiny iso}}(\sym\nabla u-\bvarepsilon_p)$\\&\\
 {\em Free energy:} &
$\frac12\,\langle\C_{\mbox{\tiny iso}}\bvarepsilon^e,\bvarepsilon^e\rangle+\,\frac12\,\mu\,
L^2_c\,\norm{\Curl(\overline{m}\,\underline{\gamma})}^2+\frac12\,\mu\, k_2\,\norm{\vecteta}^2$\\&\\
& $p=\overline{m}\,\underline{\gamma}=\sum_{\alpha}  \gamma^\alpha\,m^\alpha$\\&\\
{\em Yield condition in $\alpha$-th slip system:} &
$|\tau^\alpha_{\mbox{\tiny E}}|+g^\alpha-\yieldzero\leq0$\\
 {\em where } & $\tau^\alpha_{\mbox{\tiny E}}:=\la\Sigma_E,m^\alpha\ra$ with\\
 & $\Sigma_E:=\sigma-\mu\,L_c^2\,\Curl\Curl(\overline{m}\,\underline{\gamma})$\\
 & $g^\alpha=-\mu\,k_2\,\eta^\alpha$
 \\&\\{\em  Dissipation inequality in $\alpha$-th slip system:} &
 $\tau^\alpha_{\mbox{\tiny E}}\,\dot{\gamma}^\alpha +g^\alpha\,\dot{\eta}^\alpha\geq0$\\
 {\em Dissipation function in $\alpha$-th slip system:} &$\mathcal{D}^\alpha_{\mbox{\tiny iso}}(q^\alpha,\beta^\alpha):=\left\{\begin{array}{ll}\yieldzero |q^\alpha| &\mbox{ if } |q^\alpha|\leq\beta^\alpha,\\ \infty &\mbox{ otherwise}\end{array}\right.$\\&\\
 {\em Flow rules in primal form:} &
 $(\tau^\alpha_{\mbox{\tiny E}},g^\alpha)\in\partial \mathcal{D}^\alpha_{\mbox{\scriptsize iso}}(\dot{\gamma}^\alpha,\dot{\eta}^\alpha)$\\
{\em Flow rules in dual form:}
&$\dot{\gamma}^\alpha=\lambda^\alpha\,\dsize\frac{\tau^\alpha_{\mbox{\tiny E}}}{|\tau^\alpha_{\mbox{\tiny E}}|},\quad\qquad \dot{\eta}^\alpha=\lambda^\alpha=|\dot{\gamma}^\alpha|$\\&\\
{\em KKT conditions:} &$\lambda^\alpha\geq0$, \quad $\phi(\tau^\alpha_E,g^\alpha)\leq0$,
\quad $\lambda^\alpha\,\phi(\tau^\alpha_{\mbox{\tiny E}},g^\alpha)=0$\\&\\
 {\em Boundary conditions for $\vectgam$:} & $(\overline{m}\,\underline{\gamma})\times{n}|_{\Gamma_{\mbox{\tiny D}}}=0$,\,\, $(\Curl (\overline{m}\,\underline{\gamma}))\times{n}|_{\partial\Omega\setminus\Gamma_{\mbox{\tiny D}}}=0$\\
 {\em Function space for $\vectgam$:} & $\vectgam(t,\cdot)\in \mbox{L}^2(\Omega,\mathbb{R}^{n_{\mbox{\tiny slip}}})$,  \,\,$\Curl(\overline{m}\,\underline{\gamma})(t,\cdot)\in 
  \mbox{H}(\mbox{Curl};\;\Omega,\,\BBR^{3\times 3})$\\&\\
 \hline
\end{tabular}\caption{\footnotesize The model of single crystal gradient plasticity  with isotropic hardening.  As shown in Section \ref{wf-iso-crystal}, the weak formulation of this model is well-posed for $\underline{\gamma}(t,\cdot)\in \mbox{L}^2(\Omega,\mathbb{R}^{n_{\mbox{\tiny slip}}})$ with $\Curl(\overline{m}\,\underline{\gamma})(t,\cdot)\in 
  \mbox{L}^2(\Omega,\,\BBR^{3\times 3})$ and the boundary condition $(\overline{m}\,\underline{\gamma})\times{n}|_{\Gamma_{\mbox{\tiny D}}}=0$. The latter is well-defined in this case because $\overline{m}\,\underline{\gamma}(t,\cdot)\in 
  \mbox{H}(\mbox{Curl};\;\Omega,\,\BBR^{3\times 3})$. The condition $|q^\alpha|\leq\beta^\alpha$ in the dissipation function $\mathcal{D}^\alpha_{\mbox{\tiny iso}}$ for the $\alpha$-slip system is crucial here for the well-posednes of the model.}\label{table:isohard-crystal}\end{center}}\end{table}

\subsubsection{Weak formulation of the model}\label{wf-iso-crystal} Assume that the problem in Section \ref{strong-iso-crystal} has a solution
$(u,\vectgam,\vecteta)$. Let $v\in H^1(\Omega,\mathbb{R}^3)$
with $v_{|\Gamma_D}=0$. Multiply the equilibrium equation with
$v-\dot{u}$ and integrate in space by parts and use the
symmetry of $\sigma$ and the elasticity relation to get
\begin{equation}\label{weak-eq1-iso}
\int_{\Omega}\la\C_{\mbox{\scriptsize iso}}\sym(\nabla
u-\overline{m}\,\underline{\gamma}),\mbox{sym}(\nabla v-\nabla\dot{u})\ra\, dx=\int_\Omega
f(v-\dot{u})\,dx\, .
\end{equation}
Now,
for any $\underline{q}=(q^1,\ldots,q^{n_{\mbox{\tiny slip}}})$ with $q^\alpha\in C^\infty(\overline{\Omega})$  and any $\underline{\beta}=(\beta^1,\ldots,\beta^{n_{\mbox{\tiny slip}}})$ with $\beta^\alpha\in L^2(\Omega)$, summing (\ref{dissinequality-slip2}) over $\alpha=1,\ldots,n_{\mbox{\tiny slip}}$ and  integrating over $\Omega$, we  get
\begin{eqnarray}\label{dissinequality2-wk-iso}
\nonumber &&\hskip-1.5truecm \int_\Omega\mathcal{D}_{\mbox{\tiny iso}}(\underline{q},\underline{\beta})\,dx
-\int_\Omega\mathcal{D}_{\mbox{\tiny iso}}(\dot{\vectgam},\dot{\vecteta})\,dx -\int_\Omega
\la\C_{\mbox{\tiny iso}}\sym(\nabla u-\overline{m}\,\underline{\gamma}),\sym(\overline{m}\,\underline{q}-\overline{m}\,\dot{\underline{\gamma}})\ra\,dx\\
&&\hskip1truecm +\int_\Omega\Bigl[\mu\,k_2\,\la\vecteta,\underline{\beta}-\dot{\vecteta}\ra+\mu\,L^2_c\,\la\Curl(\overline{m}\,\underline{\gamma}),\Curl(\overline{m}\,\underline{q}-\overline{m}\,\dot{\underline{\gamma}})\ra
\,\Bigr]dx\geq0\,,\end{eqnarray}
where
\begin{equation}\label{overall-diss-iso}
\mathcal{D}_{\mbox{\tiny iso}}(\underline{q},\underline{\beta}):=\sum_\alpha \mathcal{D}^\alpha_{\mbox{\tiny iso}}(q^\alpha,\beta^\alpha)\,.\end{equation}

Now adding up (\ref{weak-eq1-iso})-(\ref{dissinequality2-wk-iso}) we
get the following weak formulation of the problem set in Section
\ref{strong-iso-crystal} in the form of a variational inequality:
\begin{eqnarray}\label{weak-form-iso-crystal}
\nonumber&&\hskip-.7truecm\int_\Omega\Bigl[\la\C_{\mbox{\tiny iso}}\sym(\nabla
u-\overline{m}\,\underline{\gamma}),\sym(\nabla v-\overline{m}\,\underline{q})-\sym(\nabla\dot{u}-\overline{m}\,\dot{\underline{\gamma}})\ra +\mu\,k_2\,\la\vecteta,\underline{\beta}-\dot{\vecteta}\ra\\
\nonumber&&\hskip6truecm+\,\, \mu\,L^2_c\,\la\Curl(\overline{m}\,\underline{\gamma}),\Curl(\overline{m}\,\underline{q}-\overline{m}\,\dot{\underline{\gamma}})\ra
\,\Bigr]dx  \\
&&\hskip2truecm +
\int_\Omega\mathcal{D}_{\mbox{\tiny iso}}(\underline{q},\underline{\beta})\,dx
-\int_\Omega\mathcal{D}_{\mbox{\tiny iso}}(\dot{\vectgam},\dot{\vecteta})\,dx\,\geq\, \int_\Omega
f\,(v-\dot{u})\,dx\,.\end{eqnarray}

\subsubsection{Existence result for the weak formulation}\label{exist-crystal-iso}
To prove the existence result for the weak formulation
(\ref{weak-form-iso-crystal}), we  follow the abstract machinery developed by
Han and Reddy in \cite{Han-ReddyBook} for mathematical problems in geometrically linear
classical plasticity and used for instance in \cite{DEMR1, REM, NCA,
EBONEFF, ENR2015} for models of gradient plasticity. To this aim, equation 
(\ref{weak-form-iso-crystal}) is written as the variational inequality of
the second kind: find $w=(u,\vectgam,\vecteta)\in \SFH^1(0,T;\SFZ)$
such that $w(0)=0$,  $\dot{w}(t)\in \SFW$ for a.e. $t\in[0,T]$ and
\begin{equation}\label{wf}
\ba(w,z-\dot{w})+j(z)-j(\dot{w})\geq \langle
\ell,z-\dot{w}\rangle\mbox{ for every } z\in \SFW\mbox{ and for a.e.
}t\in[0,T]\,,\end{equation} where $\SFZ$ is a suitable Hilbert space
and $\SFW$ is some closed, convex subset of $\SFZ$ to be constructed
later,
\begin{eqnarray}
\nonumber \ba(w,z)&=&\int_\Omega\Bigl[\la\C_{\mbox{\scriptsize iso}}\sym(\nabla
u-\overline{m}\,\underline{\gamma}),\sym(\nabla v-\overline{m}\,\underline{q})\ra+\,\mu\,k_2\,\la\vecteta,\underline{\beta}\ra\\
&&\hskip4.5truecm+ \mu\,L^2_c\,\la\Curl(\overline{m}\,\underline{\gamma}),\Curl(\overline{m}\,\underline{q})\ra
\,\Bigr]dx\,,
\label{bilin-iso-spin}
\\\nonumber\\
j(z)&=&\int_\Omega\mathcal{D}_{\mbox{\tiny iso}}(\underline{q},\underline{\beta})\,dx\,,\label{functional-isospin}\\
\langle \ell,z\rangle&=&\int_\Omega
f\,v\,dx\,,\label{lin-form}\end{eqnarray} for
$w=(u,\vectgam,\vecteta)$ and $z=(v,\underline{q},\underline{\beta})$ in
$\SFZ$.
\vskip.2truecm\noindent
The Hilbert space $\SFZ$ and the closed convex subset $\SFW$ are
constructed in such a way that the functionals $\ba$, $j$ and
$\ell$ satisfy the assumptions in the abstract result in
\cite[Theorem 6.19]{Han-ReddyBook}. The key issue here is the
coercivity of the bilinear form $\ba$ on the set $\SFW$, that is, $a(z,z)\geq C\norm{z}^2_Z$ for every $z\in\SFW$ and for some $C>0$.
\vskip.2truecm\noindent
We let
\begin{eqnarray}
\SFV&:=&\mathsf{H}^1_0(\Omega,{\Gamma_{\mbox{\scriptsize
D}}},\mathbb{R}^3)=\bigl\{v\in \mathsf{H}^1(\Omega,\mathbb{R}^3)\,\,|\,\,
v_{|\Gamma_{\mbox{\scriptsize
D}}}=0\bigr\}\,,\label{space-v-crytsal}\\
 \SFP&:=&\Bigl\{\underline{q}\in\mbox{L}^2(\Omega,\,\mathbb{R}^{n_{\mbox{\tiny slip}}})\,\,|\,\,\Curl(\overline{m}\,\underline{q})\in \mbox{L}^2(\Omega,\,\mathbb{R}^{3\times 3})\mbox{ with }(\overline{m}\,\underline{q})\times n|_{\Gamma_{\mbox{\tiny D}}}=0\Bigr\},\label{space-epsp-crystal}\\
\Lambda&:=& L^2(\Omega,\,\mathbb{R}^{n_{\mbox{\tiny slip}}})\,,\label{space-hardvar-crystal}\\
   \SFZ&:=&\SFV\times \SFP\times\Lambda\,,\label{product-space-crystal}\\
   \SFW&:=&\Bigl\{z=(v,\underline{q},\underline{\beta})\in\SFZ\,\,\,|\,\,\,|q^\alpha|\leq\beta^\alpha\,,\quad \alpha=1,\ldots,n_{\mbox{\tiny slip}}\Bigr\}\,,\label{set-W-crystal}
\end{eqnarray} and define the norms
\begin{eqnarray}
\nonumber&&  \norm{v}_V:=\norm{\nabla v}_{L^2}\,,\label{norm-V} \quad \norm{\underline{q}}^2_P:=\norm{\underline{q}}^2_{L^2} +\norm{\Curl(\overline{m}\,\underline{q})}^2_{L^2}\,,\quad  \norm{\underline{\beta}}^2_{\Lambda}=\dsize\sum_\alpha\norm{\beta^\alpha}^2_{L^2}\,,\\
&&\norm{z}^2_{Z}:=\norm{v}^2_{V} +\norm{\underline{q}}^2_{P}+
+\norm{\underline{\beta}}^2_{\Lambda}\quad\mbox{ for }
z=(v,\underline{q},\underline{\beta})\in \SFZ\,. \label{norm-Z}
\end{eqnarray}

Let us show that the bilinear form $\ba$ is coercive on $\SFW$.
 Let therefore $z=(v,\underline{q},\underline{\beta})\in \SFW$. First of all notice that 
 \begin{equation}\label{normP-dom-iso}
 \norm{\overline{m}\,\underline{q}}_{L^2}\leq\norm{\underline{q}}_P\,\leq\,\norm{\underline{\beta}}_{\Lambda}\,.\end{equation}
 So,
\begin{eqnarray}
\nonumber\ba(z,z)& \geq&  m_0\norm{\sym(\nabla v-\overline{m}\,\underline{q})}^2_{L^2}\mbox{
(from (\ref{ellipticityC}))} +\, \mu\,k_2\norm{\underline{\beta}}^2_{L^2}+\,\mu\,
L_c^2\norm{\Curl(\overline{m}\,\underline{q}}_{L^2}^2 \\
\nonumber&=&m_0\bigl[\norm{\sym\nabla v}^2_{L^2} +\norm{\sym(\overline{m}\,\underline{q})}^2_2
-2\la\sym\nabla v, \sym(\overline{m}\,\underline{q})\ra_{L^2}\bigr]\\
\nonumber &&\hskip6truecm+ \,\, \mu\,k_2\norm{\underline{\beta}}^2_{L^2}  +\,\mu L_c^2\,\norm{\Curl(\overline{m}\,\underline{q}}_{L^2}^2\\
\nonumber&\geq &m_0\left[\norm{\sym\nabla v}^2_{L^2} +\norm{\sym(\overline{m}\,\underline{q})}^2_{L^2} -\theta\norm{\sym(\nabla
v)}_{L^2}^2-\frac1\theta\norm{\sym(\overline{m}\,\underline{q})}_{L^2}^2\right]\\
\nonumber&& \hskip1truecm+\,\,\mu\,k_2\,\norm{\underline{\beta}}^2_{L^2}+\,\mu \,L_c^2\norm{\Curl(\overline{m}\,\underline{q})}_{L^2}^2\quad \mbox{ (Young's inequality for $0<\theta<1$)}\\
\nonumber &=& m_0(1-\theta)\norm{\sym\nabla
v}^2_{L^2} +m_0\Bigl(1-\frac1\theta\Bigr)\norm{\sym(\overline{m}\,\underline{q})}_{L^2}^2+\,\mu\,k_2\,\norm{\underline{\beta}}^2_{L^2}
+\mu \,L_c^2\norm{\Curl(\overline{m}\,\underline{q})}_{L^2}^2  \\
 &\geq& m_0(1-\theta)\norm{\sym\nabla
v}^2_{L^2}+\left[\frac12\,\mu\,k_2+m_0\Bigl(1-\frac1\theta\Bigr)\right]\norm{\underline{q}}_{L^2}^2     \nonumber \mbox{ \,(using both $\leq$ in  (\ref{normP-dom-iso}))}\\
 &&\hskip3truecm +\,\,\frac12\,\mu\,k_2\,\norm{\underline{\beta}}^2_{L^2}+\,\mu \,L_c^2\norm{\Curl(\overline{m}\,\underline{q})}_{L^2}^2\,.\label{coercive-a}\end{eqnarray}

Since the constant $\mu\,k_2>0$ and hence, it is possible to choose $\theta$ such that $$\displaystyle\frac{m_0}{m_0+\frac12\,\,\mu\,k_2}< \theta<1,$$  we are always able to  find some constant
$C(\theta,m_0,k,L_c,\Omega)>0$ such that
\begin{equation}\label{coerc-kin-crystal}a(z,z)\geq C\left[\norm{v}_{V^2}+\norm{\underline{q}}^2_{P} +\norm{\underline{\beta}}^2_{L^2}
\right]=C\norm{z}^2_{Z}\quad\forall z=(v,\underline{q},\underline{\beta})\in \SFW\,.\end{equation} This shows existence for the model  of single crystal gradient plasticity with isotropic hardening.
\begin{remark}\label{uniuqeness}{\rm The uniqueness result is obtained as in \cite[Section 4.4] {EHN2016} and in \cite[Section 3.2.7]{ENF2016} provided $\Curl\Curl(\overline{m}\,\vectgam)\in L^2(\Omega,\mathbb{R}^{3\times 3})$.}\end{remark}
\section{The description of the model with linear kinematic hardening}\label{single-kin}
\subsection{Constitutive equations} 
 \subsubsection{The free-energy}\label{free-eng-kin} We consider  a free
energy in the additively separated form
\begin{eqnarray}\label{free-eng}
\Psi(\nabla u,\vectgam,\Curl p):
&=&\underbrace{\Psi^{\mbox{\scriptsize lin}}_e(\sym
e)}_{\mbox{\small elastic energy}}\,\,
+\,\,\,\underbrace{\Psi^{\mbox{\scriptsize lin}}_{\mbox{\scriptsize
curl}}(\Curl p)}_{\mbox{\small defect
energy (GND)}}\\
\nonumber
&&\qquad\quad+\hskip-1.truecm\,\,\underbrace{\,\,\Psi^{\mbox{\scriptsize lin}}_{\mbox{\scriptsize
 kin}}(\vectgam)}_{\begin{array}{c}\mbox{\small kinematic}\\
\mbox{\small hardening energy (SSD)}\end{array}}\,,
 \end{eqnarray} where
 \begin{equation}\label{free-eng-expr}\left\{\begin{array}{lcl} \Psi^{\mbox{\scriptsize
lin}}_{\mbox{\scriptsize e}}(\sym e): &=&\frac12\,\la\sym e,\C_{\mbox{\scriptsize iso}} \sym
e\ra,\\\\
 \Psi^{\mbox{\scriptsize lin}}_{\mbox{\scriptsize
curl}}(\Curl p):&=&\frac12\,\mu\, L_c^2\,\norm{\Curl p}^2=\frac12\,\mu\,L_c^2\norm{\Curl(\overline{m}\,\underline{\gamma})}^2\,,\\
\\ \Psi^{\mbox{\scriptsize lin}}_{\mbox{\scriptsize
 kin}}(\vectgam)&& \mbox{ (a quadratic form to be specified later)}\,.\end{array}\right.\end{equation}

 \subsubsection{The
derivation of the dissipation inequality}\label{dis-ineq-section}
The local free-energy imbalance states that
\begin{equation}
\dot{\Psi} - \la\sigma,\dot{e}\ra - \la\sigma,\dot{p}\ra  \leq 0\
. \label{2ndlaw}
\end{equation}
Now we expand the first term, substitute (\ref{free-eng}) and get
\begin{equation}\label{exp-2ndlaw-crystal1}
\la\C_{\mbox{\scriptsize iso}}\sym
e-\sigma,\sym\dot{e}\ra-\la\sigma,\overline{m}\,\underline{\dot{\gamma}}\ra+\la\frac{\partial\Psi^{\mbox{\tiny
 lin}}_{\mbox{\tiny
 kin}}}{\partial\underline{\gamma}},\underline{\dot{\gamma}}\ra+\mu
L_c^2\la\Curl(\overline{m}\,\underline{\gamma}),\Curl(\overline{m}\,\underline{\dot{\gamma}})\ra\leq0\,,
\end{equation}
Since the inequality (\ref{exp-2ndlaw-crystal1}) must be satisfied for whatever elastic-plastic deformation mechanism, inlcuding purely elastic ones (for which  $\underline{\dot{\gamma}}=0$), inequality (\ref{exp-2ndlaw-crystal1}) implies the usual infinitesimal elastic stress-strain relation
\begin{eqnarray}
 \sigma = \C_{\mbox{\scriptsize iso}}\,\varepsilon_e\,=\,2\mu\, \sym(\nabla
u-p)+\lambda\, \tr(\nabla u-p)\id 
\,=\, 2\mu\, (\sym(\nabla
u)-\bvarepsilon_p)+\lambda\, \tr(\nabla u)\id\label{elasticlaw-kin}
\end{eqnarray}
and the local reduced dissipation inequality
\begin{equation}\label{dissip-ineq-crystal1}
\la\sigma,\overline{m}\,\underline{\dot{\gamma}}\ra-\la\frac{\partial\Psi^{\mbox{\tiny
 lin}}_{\mbox{\tiny
 kin}}}{\partial\underline{\gamma}},\underline{\dot{\gamma}}\ra-\mu
L_c^2\la\Curl(\overline{m}\,\underline{\gamma}),\Curl(\overline{m}\,\underline{\dot{\gamma}})\ra\geq0\,,
\end{equation} that we integrate over $\Omega$ and get
\begin{eqnarray}
\nonumber\hskip-4truecm0&\leq&\int_{\Omega}\Bigl[\la\sigma,\overline{m}\,\underline{\dot{\gamma}}\ra-\la\frac{\partial\Psi^{\mbox{\tiny
 lin}}_{\mbox{\tiny
 kin}}}{\partial\underline{\gamma}},\underline{\dot{\gamma}}\ra-\mu
L_c^2\la\Curl(\overline{m}\,\underline{\gamma}),\Curl(\overline{m}\,\underline{\dot{\gamma}})\ra\Bigr]\,dx\\
 \nonumber\hskip-4truecm &=&\int_{\Omega}\Bigl[\la\sigma,\overline{m}\,\underline{\dot{\gamma}}\ra-\la\frac{\partial\Psi^{\mbox{\tiny
 lin}}_{\mbox{\tiny
 kin}}}{\partial\underline{\gamma}},\underline{\dot{\gamma}}\ra-\mu
L_c^2\la\Curl\Curl(\overline{m}\,\underline{\gamma}),\overline{m}\,\underline{\dot{\gamma}}\ra\Bigr]\,dx\\
 \nonumber &&\hskip4truecm-\,\,\sum_{i=1}^3\mu\,L_c^2\int_\Omega\div\left((\overline{m}\,\underline{\dot{\gamma}})^i\times(\Curl(\overline{m}\,\underline{\gamma}))^i\right)\,dx\\
\hskip-4truecm &=&\int_{\Omega}\Bigl[\la\sigma,\overline{m}\,\underline{\dot{\gamma}}\ra-\la\frac{\partial\Psi^{\mbox{\tiny
 lin}}_{\mbox{\tiny
 kin}}}{\partial\underline{\gamma}},\underline{\dot{\gamma}}\ra-\mu
L_c^2\la\Curl\Curl(\overline{m}\,\underline{\gamma}),\overline{m}\,\underline{\dot{\gamma}}\ra\Bigr]\,dx\\
 \nonumber &&\hskip4truecm+\,\,\sum_{i=1}^3\mu\,L_c^2\int_{\partial\Omega}\la(\overline{m}\,\underline{\dot{\gamma}})^i\times(\Curl(\overline{m}\,\underline{\gamma}))^i,n\ra\,dS\,.
\end{eqnarray}
In order to obtain a dissipation
inequality in the spirit of classical plasticity, we assume that the
infinitesimal plastic distortion $p$ (and hence the slips through $\overline{m}\,\underline{\gamma}$) satisfies the so-called {\it
linearized insulation condition}
\begin{equation}\sum_{i=1}^3\mu\,L_c^2\int_{\partial\Omega}\la(\overline{m}\,\underline{\dot{\gamma}})^i\times(\Curl(\overline{m}\,\underline{\gamma}))^i,n\ra\,dS=0\,.\label{lin-sul}\end{equation} Under (\ref{lin-sul})  we then
obtain a global version of the reduced dissipation inequality
\begin{equation}\label{global-dissip-ineq1}\int_{\Omega}\Bigl[\la\sigma,\overline{m}\,\underline{\dot{\gamma}}\ra-\la\frac{\partial\Psi^{\mbox{\tiny
 lin}}_{\mbox{\tiny
 kin}}}{\partial\underline{\gamma}},\underline{\dot{\gamma}}\ra-\mu
L_c^2\la\Curl\Curl(\overline{m}\,\underline{\gamma}),\overline{m}\,\underline{\dot{\gamma}}\ra\Bigr]\,dx\geq0\,.\end{equation}
That is,
\begin{eqnarray}\label{global-dissip-ineq2}
\int_\Omega\sum_\alpha\Bigl[\la\sigma,m^\alpha\ra\,\dot{\gamma}^\alpha-\frac{\partial\Psi^{\mbox{\tiny
 lin}}_{\mbox{\tiny
 kin}}}{\partial\gamma^\alpha}\,\dot{\gamma}^\alpha-\la\mu\,L_c^2\la\Curl\Curl(\overline{m}\,\underline{\gamma}),m^\alpha\ra\,\dot{\gamma}^\alpha\Bigr]\,dx\geq0\,.\end{eqnarray} Hence, we get
 \begin{equation}\label{dissip-ineq-final}
 \int_\Omega\sum_\alpha\tau^\alpha_{\mbox{\tiny E}}\,\dot{\gamma}^\alpha\,dx\geq0\,\end{equation}
 where we set
\begin{eqnarray}\tau^\alpha_{\mbox{\tiny E}}&:=&\tau^\alpha +s^\alpha_{\mbox{\tiny loc}}+s^\alpha_{\mbox{\tiny nonloc}}\\
\mbox{ }\hskip-1truecm \tau^\alpha&:=&\la\sigma,m^\alpha\ra\quad\mbox{(resolved shear stress for the $\alpha$-th slip system)}\,,\\
s^\alpha_{\mbox{\tiny loc}}&:=&-\,\frac{\partial\Psi^{\mbox{\tiny
 lin}}_{\mbox{\tiny kin}}}{\partial\gamma^\alpha}\quad\mbox{(local backstress for the $\alpha$-th slip system)} \,,\\
s^\alpha_{\mbox{\tiny nonloc}}&:=&-\mu\,L_c^2\,\la\Curl\Curl(\overline{m}\,\underline{\gamma}),m^\alpha\ra\,,\\
\nonumber&&\mbox{(nonlocal backstress for the $\alpha$-th slip system)}
\end{eqnarray} 
 
\subsubsection{The flow rule}\label{flow-law}
We consider a yield function on the $\alpha$-th slip system
defined by
\begin{equation}
\label{yield-funct-kinematic} \phi(\tau^\alpha_{\mbox{\tiny E}}):= |\tau^\alpha_{\mbox{\tiny E}}| - \yieldzero\quad\mbox{ for }\quad\tau^\alpha_{\mbox{\tiny E}}=\tau^\alpha +s^\alpha_{\mbox{\tiny loc}}+s^\alpha_{\mbox{\tiny nonloc}}\,.
\end{equation}
Here,  $\yieldzero$ is the initial yield stress of the material, that we assume to be constant on all slip systems.
 So the set
of admissible  generalized stresses for the $\alpha$-th slip system is defined as
\begin{equation}\label{admiss-stress-kin}
\mathcal{K}^\alpha:=\left\{\tau^\alpha_{\mbox{\tiny
E}}\,\,|\,\,\phi(\tau^\alpha_{\mbox{\tiny E}})
   \leq0\right\}\,,\end{equation} with its interior $\mbox{Int}(\mathcal{K}^\alpha)$ and its boundary $\partial\mathcal{K}^\alpha$ being the generalized elastic region and the yield surface for the $\alpha$-th slip system, respectively.
   \vskip.1truecm\noindent
The  principle of maximum dissipation associated with the $\alpha$-th slip system gives  us the normality
law\begin{equation}\label{normalcone} \dot{\gamma}^\alpha\in
N_{\mathcal{K}^\alpha}(\tau^\alpha_{\mbox{\tiny E}})\,,\end{equation}
where $\dsize N_{\mathcal{K}^\alpha}(\tau^\alpha_{\mbox{\tiny E}})$
denotes the normal cone to $\mathcal{K}^\alpha$ at
$\tau^\alpha_{\mbox{\tiny E}}$. That is,  $\dot{\gamma}^\alpha$ satisfies
 \begin{equation}
(\overline{\tau}^\alpha - \tau^\alpha_{\mbox{\tiny E}})\dot{\gamma}^\alpha \leq
0\ \quad \mbox{for all}\ \overline{\tau}^\alpha \in {\mathcal{K}^\alpha}\ .
\label{normality2-alpha}
\end{equation} Notice that $N_{\mathcal{K}^\alpha}=\partial\Chi_{\mathcal{K}^\alpha}$, where
$\Chi_{\mathcal{K}^\alpha}$ denotes the indicator function of the set $\mathcal{K}^\alpha$
and $\partial\Chi_{\mathcal{K}^\alpha}$ denotes the subdifferential of the function
$\Chi_{\mathcal{K}^\alpha}$.\\ Whenever the yield surface $\partial\mathcal{K}^\alpha$ is
smooth at $\dsize\Sigma_p^\alpha$ then
$$\dot{\Gamma}^\alpha_p\in
N_{\mathcal{K}^\alpha}(\Sigma^\alpha_p)\quad\Rightarrow\quad\exists\lambda^\alpha\mbox{ such that }
\dot{\gamma}^\alpha=\lambda^\alpha\,\frac{\tau^\alpha_{\mbox{\tiny E}}}{|\tau^\alpha_{\mbox{\tiny E}}|}\quad\mbox{and hence}\quad
|\dot{\gamma}^\alpha|=\lambda^\alpha$$ with the Karush-Kuhn
Tucker conditions: $\lambda^\alpha\geq0$, $\phi(\Sigma^\alpha_p)\leq0$ and $\lambda^\alpha\,\phi(\Sigma^\alpha_p)=0$\,.\\
Using convex analysis (Legendre-transformation) we find that
\begin{eqnarray}\label{dualflowlaw-crystal} &\underbrace{\dot{\gamma}^\alpha_\in
\partial\Chi_{\mathcal{K}^\alpha}(\tau^\alpha_{\mbox{\tiny E}})}_{\mbox{\bf flow rule in its dual formulation for the $\alpha$-th slip system}}&\\
\nonumber&\Updownarrow& \\
&\underbrace{\tau^\alpha_{\mbox{\tiny E}}\in
\partial \Chi^*_{\mathcal{K}^\alpha} (\dot{\gamma}^\alpha)\,}_{\mbox{\bf flow rule in its primal formulation  for the $\alpha$-th slip system}}&\label{primalflowlaw-crystal}
\end{eqnarray} where $\Chi^*_{\mathcal{K}^\alpha}$ is the Fenchel-Legendre dual of the function $\Chi_{\mathcal{K}^\alpha}$ denoted in this context by $\mathcal{D}^\alpha_{\mbox{\scriptsize kin}}$,
 the one-homogeneous dissipation function for the $\alpha$-th slip system. That is, for every $q^\alpha$,
\begin{equation}\label{dissp-function-kin-slip}
\nonumber\mathcal{D}^\alpha_{\mbox{\scriptsize kin}}(q^\alpha)
\,=\,\sup\graffe{\tau^\alpha_{\mbox{\tiny E}}q^\alpha\,\,|\,\,\tau^\alpha_{\mbox{\tiny E}}\in\mathcal{K}^\alpha}\,=\, \yieldzero\,|q^\alpha| \,.
\end{equation} We get from the definition of the subdifferential ($\tau^\alpha_{\mbox{\tiny E}} \in
\partial \Chi^*_{\mathcal{K}^\alpha} (\dot{\gamma}^\alpha)$) that
\begin{equation}
\mathcal{D}^\alpha_{\mbox{\scriptsize kin}} (q^\alpha) \geq \mathcal{D}^\alpha_{\mbox{\scriptsize kin}}(\dot{\gamma}^\alpha) +
\tau^\alpha_{\mbox{\tiny E}}(q^\alpha-\dot{\gamma}^\alpha)\quad \mbox{for any
}q^\alpha.\label{dissineq-slip1}\end{equation} That is,
\begin{eqnarray}
\mathcal{D}^\alpha_{\mbox{\scriptsize kin}} (q^\alpha) &\geq& \mathcal{D}^\alpha_{\mbox{\scriptsize kin}}(\dot{\gamma}^\alpha) +\la\C_{\mbox{\tiny iso}}(\sym(\nabla u)-\sym(\overline{m}\,\underline{\gamma})),m^\alpha q^\alpha-m^\alpha\dot{\gamma}^\alpha\ra\quad\\
\nonumber &&\quad -\,\frac{\partial\Psi^{\mbox{\scriptsize
 lin}}_{\mbox{\tiny kin}}}{\partial\gamma^\alpha}\, (q^\alpha-\dot{\gamma}^\alpha)-\mu\,L_c^2\,\la\Curl\Curl(\overline{m}\,\underline{\gamma}),\,m^\alpha q^\alpha-m^\alpha\dot{\gamma}^\alpha \ra\qquad\mbox{for any
}q^\alpha.\label{dissineq-slip2}\end{eqnarray} 
\vskip.3truecm\noindent In the next sections, we present a complete mathematical analysis of the model including  both strong and weak formulations as well as a corresponding
existence result.
\subsection{Mathematical analysis of the model}\label{math-anal-kin}
\subsubsection{Strong formulation}\label{strong-kin-crystal} To summarize, we have obtained the following strong
formulation for the model of single crystal infinitesimal gradient plasticity with
 linear kinematic  hardening.  Given $f\in \SFH^1(0,T;L^2(\Omega,\mathbb{R}^3))$, the goal is to find:
\begin{itemize}\item[(i)] the displacement $u\in \SFH^1(0,T; H^1_0(\Omega,{\Gamma_{\mbox{\scriptsize D}}},\mathbb{R}^3))$,
\item[(ii)] the infinitesimal plastic slips $\gamma^\alpha\in
\SFH^1(0,T;L^2(\Omega))$ for $\alpha=1,\ldots,n_{\mbox{\tiny slip}}$ with \\ $\Curl(\overline{m}\,\underline{\gamma})\in\SFH^1(0,T;H(\Curl;\Omega, \mathbb{R}^{3\times 3})$  \end{itemize}
 such that the content of Table \ref{table:micro-isohard-crystal} holds.
 
 \begin{table}[ht!]{\footnotesize\begin{center}
\begin{tabular}{|ll|}\hline &\qquad\\
 {\em Additive split of distortion:}& $\nabla u =e +p$,\quad $\bvarepsilon_e=\mbox{sym}\,e$,\quad $\bvarepsilon_p=\sym p$\\
{\em Plastic distortion in slip system:}&$p=\dsize\sum_{\alpha=1}^{n_{\mbox{\tiny silp}}}      \gamma^\alpha\,m^\alpha$ with $m^\alpha= l^\alpha\otimes \nu^\alpha$, \quad $\tr(p)=0$\\
{\em Equilibrium:} & $\mbox{Div}\,\sigma +f=0$ with
$\sigma=\C_{\mbox{\tiny iso}}\bvarepsilon_e=\C_{\mbox{\tiny iso}}(\sym\nabla u-\bvarepsilon_p)$
\\&\\
 {\em Free energy:} &
$\frac12\,\langle\C_{\mbox{\tiny iso}}\bvarepsilon^e,\bvarepsilon^e\rangle+\Psi^{\mbox{\tiny lin}}_{\mbox{\tiny
 kin}}(\vectgam) +\,\frac12\,\mu\,
L^2_c\,\norm{\Curl (\overline{m}\,\vectgam)}^2$\\&\\
& $ \Psi^{\mbox{\tiny
 lin}}_{\mbox{\tiny
 kin}}(\vectgam)$ to be specified later; \,\,\, $p=\overline{m}\,\underline{\gamma}=\sum_{\alpha}  \gamma^\alpha\,m^\alpha$\\&\\
{\em Yield condition in $\alpha$-th slip system:} &
$|\tau^\alpha_{\mbox{\tiny E}}|-\yieldzero\leq0$\\
 {\em where } & $\tau^\alpha_{\mbox{\tiny E}}:=\la\sigma -\mu\,L_c^2\,\Curl\Curl(\overline{m}\,\underline{\gamma}),m^\alpha\ra +\dsize\frac{\partial\Psi^{\mbox{\tiny
 lin}}_{\mbox{\tiny kin}}}{\partial\gamma^\alpha}$ 
 \\&\\{\em  Dissipation inequality in $\alpha$-th slip system:} &
 $\tau^\alpha_{\mbox{\tiny E}}\,\dot{\gamma}^\alpha\geq0$\\
 {\em Dissipation function in $\alpha$-th slip system:} &$\mathcal{D}^\alpha_{\mbox{\tiny kin}}(q^\alpha):=\yieldzero |q^\alpha|$\\&\\
 {\em Flow rules in primal form:} &
 $\tau^\alpha_{\mbox{\tiny E}}\in\partial \mathcal{D}^\alpha_{\mbox{\tiny kin}}(\dot{\gamma}^\alpha)$\\&\\
{\em Flow rules in dual form:}
&$\dot{\gamma}^\alpha=\lambda^\alpha\,\dsize\frac{\tau^\alpha_{\mbox{\tiny E}}}{|\tau^\alpha_{\mbox{\tiny E}}|},\quad\qquad \lambda^\alpha=|\dot{\gamma}^\alpha|$\\&\\
{\em KKT conditions:} &$\lambda^\alpha\geq0$, \quad $\phi(\tau^\alpha_E)\leq0$,
\quad $\lambda^\alpha\,\phi(\tau^\alpha_{\mbox{\tiny E}})=0$\\&\\
 {\em Boundary conditions for $\underline{\gamma}$:} & $(\overline{m}\,\underline{\gamma})\times{n}|_{\Gamma_{\mbox{\tiny D}}}=0$,\,\, $(\Curl (\overline{m}\,\underline{\gamma}))\times{n}|_{\partial\Omega\setminus\Gamma_{\mbox{\tiny D}}}=0$\\&\\
 {\em Function space for $\underline{\gamma}$:} & $\vectgam(t,\cdot)\in L^2(\Omega,\mathbb{R}^{n_{\mbox{\tiny slip}}}),\,\quad (\overline{m}\,\underline{\gamma})(t,\cdot)\in \mbox{H}(\mbox{Curl};\;\Omega,\,\BBR^{3\times 3})$\\&\\
 \hline
\end{tabular}\caption{\footnotesize The model  of single crystal gradient plasticity  with linear kinematic hardening. As shown in Section \ref{wf-kin-crystal}, the weak formulation of this model is also well-posed for $\underline{\gamma}(t,\cdot)\in \mbox{L}^2(\Omega,\mathbb{R}^{n_{\mbox{\tiny slip}}})$ with $\Curl(\overline{m}\,\underline{\gamma})(t,\cdot)\in 
  \mbox{L}^2(\Omega,\,\BBR^{3\times 3})$ and the boundary condition $(\overline{m}\,\underline{\gamma})\times{n}|_{\Gamma_{\mbox{\tiny D}}}=0$.}\label{table:micro-isohard-crystal}\end{center}}\end{table}

\subsubsection{Weak formulation}\label{wf-kin-crystal} Assume that the problem in Section \ref{strong-kin-crystal} has a solution
$(u,\vectgam,\vecteta)$.  Let $v\in H^1(\Omega,\mathbb{R}^3)$
with $v_{|\Gamma_D}=0$. Multiply the equilibrium equation with
$v-\dot{u}$ and integrate in space by parts and use the
symmetry of $\sigma$ and the elasticity relation to get
\begin{equation}\label{weak-eq1}
\int_{\Omega}\la\C_{\mbox{\scriptsize iso}}\sym(\nabla
u-\overline{m}\,\underline{\gamma}),\mbox{sym}(\nabla v-\nabla\dot{u})\ra\, dx=\int_\Omega
f(v-\dot{u})\,dx\, .
\end{equation}
Now, for any $\underline{q}=(q^1,\ldots,q^{n_{\mbox{\tiny slip}}})$ with $q^\alpha\in C^\infty(\overline{\Omega})$, summing (\ref{dissineq-slip2}) over $\alpha=1,\ldots,n_{\mbox{\tiny slip}}$, then integrating over $\Omega$, in particular integrating by
parts the term with Curl\,Curl using the boundary conditions
$$(\overline{m}\,\underline{q}-\overline{m}\,\underline{\dot{\gamma}})\times\,n=0\mbox{ on }\Gamma_D,\qquad
\mbox{Curl}(\overline{m}\,\underline{\gamma})\times\,n=0\mbox{ on }
\partial\Omega\setminus\Gamma_D$$ and get, 
\begin{eqnarray}\label{dissinequality2-wk}
\nonumber &&\hskip-1truecm \int_\Omega\mathcal{D}_{\mbox{\tiny kin}}(\underline{q})\,dx
-\int_\Omega\mathcal{D}_{\mbox{\tiny kin}}(\dot{\vectgam})\,dx -\int_\Omega
\la\C_{\mbox{\tiny iso}}\sym(\nabla u-\overline{m}\,\underline{\gamma}),\sym(\overline{m}\,\underline{q}-\overline{m}\,\dot{\underline{\gamma}})\ra\,dx\\
&&\quad-\int_\Omega\Bigl[\la\frac{\partial\Psi_{\mbox{\tiny kin}}}{\partial\underline{\gamma}}, \underline{q}-\underline{\dot{\gamma}}\ra+\mu\,L_c^2\,\la\Curl(\overline{m}\,\underline{\gamma}),\,\Curl(\underline{m}\,\underline{q}-\overline{m}\,\underline{\dot{\gamma}}) \ra\Bigr]\,dx
\end{eqnarray}
where
\begin{equation}\label{overall-diss}
\mathcal{D}_{\mbox{\tiny iso}}(\underline{q}):=\sum_\alpha \mathcal{D}^\alpha_{\mbox{\tiny iso}}(q^\alpha)\,.\end{equation}

Now adding up (\ref{weak-eq1}) and (\ref{dissinequality2-wk}) we
get the following weak formulation of the problem in Section
\ref{strong-kin-crystal} in the form of a variational inequality:
\begin{eqnarray}\label{weak-form-kin}
\nonumber&&\hskip-1truecm\int_\Omega\Bigl[\la\C_{\mbox{\scriptsize iso}}(\sym(\nabla
u-\overline{m}\,\underline{\gamma})),\sym(\nabla v-\overline{m}\,\underline{q})-\sym(\nabla\dot{u}-\overline{m}\,\dot{\underline{\gamma}})\ra\,\\
&&\hskip2truecm +\,\,\la\frac{\partial\Psi_{\mbox{\tiny kin}}}{\partial\underline{\gamma}}, \underline{q}-\underline{\dot{\gamma}}\ra+\,\,\mu\,L_c^2\la\Curl(\overline{m}\,\underline{\gamma}),\Curl(\overline{m}\,\underline{q}-\overline{m}\,\dot{\underline{\gamma}})\ra\Bigr]\,dx\\
&&\hskip4truecm +\int_\Omega\mathcal{D}_{\mbox{\tiny kin}}(\underline{q})\,dx
-\int_\Omega\mathcal{D}_{\mbox{\tiny kin}}(\dot{\vectgam})\,dx\,\geq\,\int_\Omega
f(v-\dot{u})\,dx\,.\nonumber\end{eqnarray}
\subsubsection{Existence result for the weak formulation}\label{existence}
 As in Section \ref{exist-crystal-iso}, we write 
(\ref{weak-form-kin})  as a variational inequality of
the second kind: find $w=(u,\vectgam)\in \SFH^1(0,T;\SFZ)$
such that $w(0)=0$  and
\begin{equation}\label{wf}
\ba(\dot{w},z-w)+j_0(z)-j_0(\dot{w})\geq \langle
\ell,z-\dot{w}\rangle\mbox{ for every } z\in \SFZ\mbox{ and for a.e.
}t\in[0,T]\,,\end{equation} where $\SFZ$ is a suitable Hilbert space
 to be constructed
later,

\begin{eqnarray}
\nonumber \ba(w,z)&=&\int_\Omega\Bigl[\la\C_{\mbox{\scriptsize iso}}\sym(\nabla
u-\overline{m}\,\underline{\gamma}),\sym(\nabla v-\overline{m}\,\underline{q})\ra+\la\frac{\partial\Psi_{\mbox{\tiny kin}}}{\partial\underline{\gamma}}, \underline{q}\ra\\
&&\hskip4truecm +\,\,\mu\,L_c^2\la\Curl(\overline{m}\,\underline{\gamma}),\Curl(\overline{m}\,\underline{q})\ra\Bigr]\,dx\,,
\label{bilin-iso-spin}
\\\nonumber\\
j(z)&=&\int_\Omega\mathcal{D}_{\mbox{\tiny kin}}(\underline{q})\,dx\,,\label{functional-isospin}\\
\langle \ell,z\rangle&=&\int_\Omega
f\,v\,dx\,,\label{lin-form}\end{eqnarray} for
$w=(u,\vectgam)$ and $z=(v,\underline{q})$ in
$\SFZ$.
\vskip.2truecm\noindent
The Hilbert space $\SFZ$ is 
constructed in such a way that the functionals $\ba$, $j$ and
$\ell$ satisfy the assumptions in the abstract result in
\cite[Theorem 6.19]{Han-ReddyBook}. The key issue here is the
coercivity of the bilinear form $\ba$ on the space $\SFZ$, that is, $a(z,z)\geq C\norm{z}^2_Z$ for every $z\in\SFZ$ and for some $C>0$. \\
While the function space for the displacement field is the same as in the polycrystalline setting (\cite{DEMR1, REM, NCA,
EBONEFF, ENR2015,EHN2016}), that is,
$$\SFV:=\mathsf{H}^1_0(\Omega,{\Gamma_{\mbox{\scriptsize
D}}},\mathbb{R}^3)=\bigl\{v\in \mathsf{H}^1(\Omega,\mathbb{R}^3)\,\,|\,\,
v_{|\Gamma_{\mbox{\scriptsize
D}}}=0\bigr\}\,,$$
the choice of the space of plastic slips requires some assumptions on the kinematic hardening density.
\begin{remark}\label{no-usual-kin}{\rm
Notice that the usual Prager-type linear kinematic hardening energy density from isotropic polycrystalline plasticity 
\begin{equation}\label{kin-poly} \Psi^{\mbox{\tiny lin}}_{\mbox{\tiny kin}}(\vectgam)=\frac12\,\mu\,k_1\norm{\sym(\overline{m}\,\vectgam)}^2=\frac12\,\mu\,k_1\,\norm{\bvarepsilon_p}^2\end{equation} and the usual  
$\mbox{H}_0(\mbox{Curl};\,\Omega,\,\Gamma_{\mbox{\scriptsize
D}},\,\sL(3))$-space (defined in (\ref{spacep-bc})) for the plastic distortion $p=\overline{m}\,\vectgam$ and its norm
\begin{equation}\label{normp-poly}\norm{\overline{m}\,\vectgam}^2_{\mbox{\scriptsize H}(\mbox{\scriptsize
Curl};\Omega)}=\norm{\overline{m}\,\vectgam}^2_{L^2}+\norm{\Curl(\overline{m}\,\vectgam)}^2_{L^2}\end{equation} are not appropriate in this context because  the set of plastic distortions
$$\Bigl\{p=\sum_\alpha \gamma^\alpha\,m^\alpha\Bigr\}$$ is not closed in the $H(\Curl)$-norm.  This means that, given a sequence $ (p_n)_n$ which converges to some $p$ in the norm $\norm{\cdot}_{\mbox{\scriptsize H}(\mbox{\scriptsize
Curl};\Omega)}$ with $ p_n=\sum_\alpha \gamma_n^\alpha\,m^\alpha$ there is no reason why $p=\sum_\alpha \gamma^\alpha\,m^\alpha$  for some $\underline{\gamma}$. \\
Also, (\ref{normp-poly}) does not define a norm in the slips $\gamma^\alpha$, $\alpha=1,\ldots,n_{\mbox{\tiny slip}}$.}\end{remark}
The natural condition on the kinematic hardening enegry density for the mathematical well-posedness of the model is $\Psi_{\mbox{\tiny kin}}$ to be quadratic and positive definite in $\vectgam$. This condition implies that
\begin{equation}\label{choice-psi-kin}
\Psi^{\mbox{\tiny lin}}_{\mbox{\tiny kin}}(\vectgam)\geq k\,\norm{\vectgam}^2= k\,\sum_{\alpha}|\gamma^\alpha|^2\,\quad\mbox{ for some constant }k>0\,\end{equation}
and therefore the choice of the space of slips will be
\begin{equation}\label{slip-space}
\SFP:=\Bigl\{\underline{q}=(q^1,\ldots,q^{n_{\mbox{\tiny slip}}})\in  L^2(\Omega,\,\mathbb{R}^{n_{\mbox{\tiny slip}}})\,\,\bigl|\,\,\mbox{Curl\,}(\overline{m}\,\underline{q})\in
L^2(\Omega,\,\BBR^{3\times 3})\Bigr\}\,\end{equation} equipped with the norm
\begin{equation}\label{norm-slip-space} 
\norm{\underline{q}}_P^2:=\norm{\underline{q}}^2_{L^2(\Omega,\mathbb{R}^{n_{\mbox{\tiny slip}}})}+\norm{\Curl(\overline{m}\,\underline{q})}^2_{L^2(\Omega,\mathbb{R}^{3\times 3})}\,.
\end{equation}
We set \begin{eqnarray}
\SFZ:&=&\SFV\times\SFP\\
\norm{z}^2_Z:&=&\norm{v}^2_V+\norm{\underline{q}}^2_P\qquad\forall z=(v,\underline{q})\in\SFZ\,.\end{eqnarray}
Now, let u show that the bilinear form $\ba$ is $\SFZ$-coercive. Let $z=(v,\underline{q})\in Z$.  First of all notice that 
 \begin{equation}\label{normP-dom}
 \norm{\sym(\overline{m}\,\underline{q})}_{L^2}\leq \norm{\overline{m}\,\underline{q}}_{L^2}\leq\norm{\underline{q}}_P\,.\end{equation}
 So,
\begin{eqnarray}
\nonumber\ba(z,z)& \geq&  m_0\norm{\sym(\nabla v-\overline{m}\,\underline{q})}^2_{L^2}\mbox{
(from (\ref{ellipticityC}))} +\, \la\frac{\partial\Psi_{\mbox{\tiny kin}}}{\partial\underline{q}}(\underline{q}), \underline{q}\ra+\,\mu\,
L_c^2\norm{\Curl(\overline{m}\,\underline{q}}_{L^2}^2 \\
\nonumber&=&m_0\bigl[\norm{\sym\nabla v}^2_{L^2} +\norm{\sym(\overline{m}\,\underline{q})}^2_2
-2\la\sym\nabla v, \sym(\overline{m}\,\underline{q})\ra_{L^2}\bigr]\\
\nonumber &&\hskip6truecm+ \,\,\la\frac{\partial\Psi_{\mbox{\tiny kin}}}{\partial\underline{q}}(\underline{q}), \underline{q}\ra  +\,\mu L_c^2\,\norm{\Curl(\overline{m}\,\underline{q}}_{L^2}^2\\
\nonumber&\geq &m_0\left[\norm{\sym\nabla v}^2_{L^2} +\norm{\sym(\overline{m}\,\underline{q})}^2_{L^2} -\theta\norm{\sym(\nabla
v)}_{L^2}^2-\frac1\theta\norm{\sym(\overline{m}\,\underline{q})}_{L^2}^2\right]\\
\nonumber&& \hskip1truecm+\,\,k\,\norm{\underline{q}}^2_{L^2}+\,\mu \,L_c^2\norm{\Curl(\overline{m}\,\underline{q})}_{L^2}^2\quad \mbox{ (Young's inequality for $0<\theta<1$)}\\
\nonumber &=& m_0(1-\theta)\norm{\sym\nabla
v}^2_{L^2} +m_0\Bigl(1-\frac1\theta\Bigr)\norm{\sym(\overline{m}\,\underline{q})}_{L^2}^2+\,k\,\norm{\underline{q}}^2_{L^2}
+\mu \,L_c^2\norm{\Curl(\overline{m}\,\underline{q})}_{L^2}^2  \\
 &\geq& m_0(1-\theta)\norm{\sym\nabla
v}^2_{L^2}+\left[k+m_0\Bigl(1-\frac1\theta\Bigr)\right]\norm{\underline{q}}_{L^2}^2     \nonumber \mbox{ \,(using (\ref{normP-dom}))} +\,\mu \,L_c^2\norm{\Curl(\overline{m}\,\underline{q})}_{L^2}^2\,.\label{coercive-a}\end{eqnarray}

Since the constant $k>0$ and hence, it is possible to choose $\theta$ such that $$\displaystyle\frac{m_0}{m_0+k}< \theta<1,$$  we are always able to  find some constant
$C(\theta,m_0,k,L_c,\Omega)>0$ such that
\begin{equation}\label{coerc-kin-crystal}a(z,z)\geq C\left[\norm{v}_{V^2}+\norm{\underline{q}}^2_{P}
\right]=C\norm{z}^2_{Z}\quad\forall z=(v,\underline{q})\in \SFZ\,.\end{equation} This shows existence {\bf and} unqiueness  for the model  of single crystal gradient plasticity with linear kinematical hardening.
\begin{remark}\label{rem-choice-psi-kin}{\rm {\bf (Some examples of $\Psi^{\mbox{\scriptsize
 lin}}_{\mbox{\tiny kin}}$ satisfying (\ref{choice-psi-kin}))}\begin{itemize}\item[(i)] Notice that even in the case where the slip planes are mutually orthogonal and this implies that $n_{\mbox{\tiny slip}}=3$, the classical Prager-type linear kinematic energy density might still not sastisfy the inequality (\ref{choice-psi-kin}).  In fact, in that case we obtain
 \begin{eqnarray}\label{choice-psi-kin-eq1}
 \nonumber \norm{\sym(\overline{m}\,\underline{\gamma})}^2&=&\la\sym(\overline{m}\,\underline{\gamma}),\sym(\overline{m}\,\underline{\gamma})\ra\,=\,\la\sym(\overline{m}\,\underline{\gamma}),\overline{m}\,\underline{\gamma}\ra\\
 \nonumber&=&\frac12\sum_{\alpha,\beta}\la\gamma^\alpha( l^\alpha\otimes\nu^\alpha+\nu^\alpha\otimes l^\alpha),\gamma^\beta\,l^\beta\otimes\nu^\beta\ra\\
 \nonumber &=&\frac12\sum_{\alpha,\beta}\gamma^\alpha\,\gamma^\beta\Bigl[\la l^\alpha,l^\beta\ra\,\la\nu^\alpha,\nu^\beta\ra +\la l^\alpha,\nu^\beta\ra\,\la \nu^\alpha,l^\beta\ra\Bigr]\\
&=& \frac12\sum_\alpha|\gamma^\alpha|^2+\underbrace{\sum_{\alpha<\beta}\gamma^\alpha\,\gamma^\beta\la l^\alpha,\nu^\beta\ra\,\la \nu^\alpha,l^\beta\ra}_{\mbox{with an undetermined sign}}\,.
 \end{eqnarray}
 However, if we rather assume that the slip systems $\{ l^\alpha,\,\nu^\alpha\}_{\alpha=1}^3$ are mutually orthogonal, which implies that the slip planes are mutually orthogonal and that $l^\alpha\perp\nu^\beta$, then in that case the expression with an undetermined sign in (\ref{choice-psi-kin-eq1}) vanishes so that we get
 \begin{equation}\label{ortho-systems}\norm{\sym(\overline{m}\,\underline{\gamma})}^2=\frac12\,\sum_\alpha|\gamma^\alpha|^2
\end{equation} and  (\ref{choice-psi-kin}) is satisfied. This situation corresponds for instance to the three coordinate planes in $\mathbb{R}^3$ with the slip directions $\{l^\alpha\}_\alpha$ suitably chosen along the three axes.
\item[(ii)] In general, we can consider 
\begin{equation}\label{kin-gen} \Psi^{\mbox{\tiny lin}}_{\mbox{\tiny kin}}(\vectgam)=\frac12\,\la\mathbb{H}\,\vectgam,\vectgam\ra
\end{equation} where $\mathbb{H}$ is a positive definite matrix in $\mathbb{R}^{n_{\mbox{\tiny slip}}\times n_{\mbox{\tiny slip}}}$.

 \end{itemize}

}\end{remark}
\section{Open problems}\label{conclusion}
 In the paper \cite{ENF2016}, we have treated a formal micromorphic penalty regularization of the models presented here. These models are considerably simpler. It would be interesting to make the limit passage rigorous, with attending convergence estimates. Since the ultimimate goal is the passage from single to polycrystalline samples, we face the following problem already alluded to in the introduction: any standard single crystal formulation (also for strain gradient plasticity) imparts a strict control on the level of the glide systems, while the polycrystalline setting must allow for some freedom w.r.t. infinitesimal plastic  rotations, while {\bf not} suppressing them. In order to reconcile these two aspects in certain respects, it is tempting to use a linear kinematic hardening contribution 
 \begin{equation}\label{lin-kin-contrib} \frac12\,\mu\,k_1\,\norm{\sym(\overline{m}\,\vectgam)}^2\,,\end{equation}
in the single crystal case. As we saw in Remark \ref{rem-choice-psi-kin}, this does not yet fit into our presented mathematical framework. However, our belief is that the variant (\ref{lin-kin-contrib}) should work as well. This seems to need an entirely new approach, perhaps based on the new Korn's inequalities for incompatible  tensor fields established  in \cite{NPW2011-1, NPW2012-1, NPW2012-2,NPW2014}.

\end{document}